\documentclass[12pt,reqno]{amsart}
\usepackage{amssymb,amsthm,amsfonts,amsmath,amscd,amssymb,epsfig,psfrag}
\usepackage{calrsfs}
\usepackage[all]{xy}

\usepackage{enumerate,array, graphicx,color,amsthm,
bbm,mathrsfs}                      
\usepackage{a4wide}           
\theoremstyle{plain}
\newtheorem{theorem}{Theorem}[section]
\newtheorem{lemma}[theorem]{Lemma}
\newtheorem{proposition}[theorem]{Proposition}
\newtheorem{corollary}[theorem]{Corollary}
\newtheorem*{remark*}{Remark}
\newtheorem*{remarks*}{Remarks}
\newtheorem{remark}[theorem]{Remark}
\newtheorem{remarks}[theorem]{Remarks}

\newtheorem*{example*}{Example}                           
\newtheorem*{examples*}{Examples}
\newtheorem*{examplecon*}{Example (continued)}                
\newtheorem{definition}[theorem]{Definition}

\DeclareMathAlphabet{\pazocal}{OMS}{zplm}{m}{n}

\newcommand{\eps}{\varepsilon}

\newcommand{\NN}{{\mathbb{N}}}
\newcommand{\ZZ}{{\mathbb{Z}}}
\newcommand{\RR}{{\mathbb{R}}}
\newcommand{\CC}{{\mathbb{C}}}

\newcommand{\TT}{{\mathbb{T}}}

\newcommand{\m}{{\bf m}}

\newcommand{\ind}{{\rm ind}}
\newcommand{\coker}{{\rm coker }}  
\newcommand{\im}{{\rm im }}        

\renewcommand{\max}{{\rm max}}

\newcommand{\Hess}{{\rm Hess}}
\newcommand{\Crit}{{\rm Crit}\,}

\newcommand{\con}{{\rm con}}
\newcommand{\reg}{{\rm reg}}

\newcommand{\loc}{{\rm loc}}
\newcommand{\ev}{{\rm ev}}

\newcommand{\FH}{{\rm FH}}
\newcommand{\FC}{{\rm FC}}
\newcommand{\MH}{{\rm MH}}
 
\newcommand{\RFH}{{\rm RFH}}
\newcommand{\SH}{{\rm SH}}
\newcommand{\supp}{{\rm supp \2}}
\newcommand{\spa}{{\rm span}}

\newcommand{\diag}{{\rm diag}}

\def\cA{{\pazocal A}}

\def\cE{{\pazocal E}}
\def\cF{{\pazocal F}}
\def\cG{{\pazocal G}}
\def\cH{{\pazocal H}}
\def\cJ{{\pazocal J}}
\def\cL{{\pazocal L}}
\def\cM{{\pazocal M}}
\def\cN{{\pazocal N}}

\def\cP{{\pazocal P}}

\def\cV{{\pazocal V}}

\def\bJ{{\mathbf J}}

\newcommand{\proofend}{\hspace*{\fill} $\Box$\\}
\newcommand{\diam}{\hspace*{\fill} $\Diamond$}

\def\1{\:\!}
\def\2{\;\!}
\def\ni{\noindent}
\def\b{\bigskip}
\def\m{\medskip}
\def\pp{\partial}
\def\id{\operatorname{id}}
\def\sma{\operatorname{small}}
\def\lar{\operatorname{large}}

\begin{document}

\title[$S^1$-equivariant Rabinowitz--Floer homology]
{$S^1$-equivariant Rabinowitz--Floer homology}

\author{Urs Frauenfelder}
\address{Urs Frauenfelder,
    Mathematisches Institut, Universit\"at Augsburg}
\email{urs.frauenfelder@math.uni-augsburg.de}

\author{Felix Schlenk}  
\thanks{FS partially supported by SNF grant 200020-144432/1.}
\address{Felix Schlenk,
Institut de Math\'ematiques,
Universit\'e de Neuch\^atel}
\email{schlenk@unine.ch}


\date{\today}

\begin{abstract}
We define the $S^1$-equivariant Rabinowitz--Floer homology of a bounding 
contact hypersurface~$\Sigma$ in an exact symplectic manifold,  
and show by a geometric argument that it vanishes if~$\Sigma$ is displaceable.

In the appendix we describe an approach to transversality for Floer homologies
for which the moduli space~$\widehat \cM_J$ of all gradient flow lines is compact for some almost complex structure~$J$. 
This approach uses a large set of perturbations, namely vector fields on the loop space, 
and selects from the possibly non-compact perturbed moduli spaces
a part near~$\widehat \cM_J$ that turns out to be compact for small enough perturbations.
\end{abstract}

\maketitle
\tableofcontents


\section{Introduction}  \label{s:intro}

Consider a bounding contact hypersurface $\Sigma$ in an exact convex symplectic manifold $(M,\lambda)$.
(Definitions are recalled in Section~\ref{s2}.)
In this situation, Kai Cieliebak and the first author defined in~\cite{CieFra09} 
a homology group $\RFH (\Sigma,M)$, the Rabinowitz--Floer homology of~$\Sigma$,
as the Floer homology associated to the Rabinowitz action functional
$$
\cA^F \colon \cL \times \RR \to \RR, \qquad
(v,\eta) \,\mapsto\, -\int_{S^1} v^*\lambda - \eta \int_{S^1} F \bigl( v(t) \bigr) \,dt .
$$
Here, $F \colon M \to \RR$ is a suitable function with $F^{-1}(0)=\Sigma$,
and $S^1 = \RR/\ZZ$ denotes the circle and 
$\cL = C^\infty(S^1,M)$ the free loop space of~$M$. 
Note that the Rabinowitz action functional is invariant under the circle action
$\tau v(\cdot) \mapsto v(\cdot -\tau)$ obtained by rotating the loop~$v$.
This makes it possible to construct the equivariant Rabinowitz--Floer homology 
$\RFH^{S^1} (\Sigma, M)$ as well.

Recall that $\Sigma$ is said to be Hamiltonian displaceable if there exists a
compactly supported Hamiltonian diffeomorphism that disjoins~$\Sigma$ from itself. 
One of the most useful properties of the Rabinowitz--Floer homology of~$\Sigma$ is that it vanishes 
if~$\Sigma$ is displaceable. 
The main result of this note is that this fact continues to hold in the equivariant case.

\m \ni
\textbf{Theorem A.}\, 
\emph{Assume that $\Sigma$ is Hamiltonian displaceable. 
Then \,$\RFH^{S^1}(\Sigma,M) = \{0\}$.}

\m
We shall prove this result by a leafwise intersection argument, 
following~\cite{AlbFra10}.
A more algebraic proof of Theorem~A was given in~\cite{BouOan12} in the framework
of symplectic homology, and their proof should also apply to Rabinowitz--Floer homology, 
cf.~Section~\ref{s:other}.

\b 
The main body of this note is organized as follows.
In Section~\ref{s2} we recall the construction of non-equivariant Rabinowitz--Floer homology $\RFH(\Sigma,M)$,
and in Section~\ref{s3} we construct $S^1$-equivariant Rabinowitz--Floer homology $\RFH^{S^1}(\Sigma,M)$.
The core of this part is Section~\ref{s:proofThA} in which we prove 
Theorem~A.
In Section~\ref{s:invariance} we give an alternative and somewhat easier approach to the invariance 
of~$\RFH^{S^1}(\Sigma,M)$.
In Section~\ref{s:other} we briefly discuss other approaches to proving 
$\RFH^{S^1}(\Sigma,M)=0$ for displaceable hypersurfaces. 
 
\b
\ni
{\bf Transversality.}
Consider a Floer theory on a symplectic manifold $(M,\omega)$ with action functional~$\cA$, 
fix an $\omega$-compatible almost complex structure~$J$ on~$M$, 
and let $\nabla \cA$ be the $L^2$-gradient of~$\cA$ defined by~$J$.
Fix $a<b$ and let $\Crit \cA_a^b$ be the critical points of~$\cA$ with action in $[a,b]$.
Denote by $\mathcal{G}_a^b$ the space of all flow lines of~$-\nabla \cA$ between elements in~$\Crit \cA_a^b$,
and by $\mathcal{G}(c_-,c_+)$ its subset of flow lines between two critical points $c_-,c_+$.
For a perturbation~$v$ of~$-\nabla \cA$ that vanishes at $\Crit \cA_a^b$, 
denote by $\mathcal{G}_a^b(v)$ the space of all flow lines of~$-\nabla \cA+v$ between elements in~$\Crit \cA_a^b$,
and by $\mathcal{G}(c_-,c_+,v)$ its subset of flow lines between $c_-,c_+$.
The transversality problem in this situation is to show that for a generic choice of perturbations~$v$, 
the spaces~$\mathcal{G}(c_-,c_+,v)$ are cut out transversally from 
``the space of all lines'' from $c_-$ to~$c_+$, and are thus smooth manifolds.

\subsubsection*{Geometric perturbations.}
The traditional approach to achieve transversality is to choose a very small space of perturbations,
namely the space of $\omega$-compatible almost complex structures on~$M$,
that may also depend on $t \in S^1$ (and for~$\RFH$ also on~$\eta \in \RR$), see~\cite{FHS:95, BouOan10}.
The advantage of this approach is that if compactness of~$\mathcal{G}_a^b$ holds, then it holds
for all spaces $\mathcal{G}_a^b(v)$ for the same reason;
its disadvantage is that one has few perturbations at hand, which makes it sometimes very hard to achieve 
transversality.
This approach has been carried out for Rabinowitz--Floer homology in~\cite[\S~4]{AbbMer15},
and the additional compact $S^{2N+1}$-factor of the $S^1$-equivariant theory causes no additional problems.

\subsubsection*{Polyfolds.}
A more conceptual approach to transversality, that applies in very general situations, 
is the theory of polyfolds~\cite{Hof06,HWZ14}. The group~$\mathbb{R}$ acts on the spaces
$\mathcal{G}(c_-,c_+)$ by time shift, and the quotient $\mathcal{G}(c_-,c_+)/\mathbb{R}$ of unparametrized gradient flow lines
can be compactified to $\overline{\mathcal{G}(c_-,c_+)/\mathbb{R}}$ by adding the unparametrized broken gradient flow lines from~$c_-$ to~$c_+$. 
While the space $\mathcal{G}(c_-,c_+)$ can be interpreted as the zero-set of a section from a Hilbert manifold
to a Hilbert bundle, this is not possible anymore for $\overline{\mathcal{G}(c_-,c_+)/\mathbb{R}}$
for two reasons: First of all the $\RR$-action on $\mathcal{G}(c_-,c_+)$ is not smooth 
in the usual sense. However, it is smooth in a new sense discovered by Hofer, Wysocki and Zehnder, namely scale smooth. 
Scale smoothness does not require just one Hilbert manifold, but a whole scale of Hilbert manifolds, 
and therefore leads to exciting interactions between Floer homology and interpolation theory~\cite{Tri78}. 
The second issue is the presence of broken gradient flow lines, which is an analytical limit phenomenon. 
However, the space of all unparametrized broken or unbroken lines from $c_-$ to~$c_+$ 
can be interpreted as the fixed point set 
of a scale smooth retraction. 
(This is in sharp contrast to the Hilbert set-up, where by the last theorem of Cartan the fixed point set of 
a smooth retraction is itself a Hilbert manifold~\cite{Car86}.) 
These two facts allow to interpret the moduli space $\overline{\mathcal{G}(c_-,c_+)/\mathbb{R}}$
as the zero-set of a Fredholm section from an M-polyfold to an M-polyfold bundle. 
Here, M stands for ``manifold flavoured", indicating that no orbifold technology is required. 
A detailed account of this story is currently written up by Albers and Wysocki~\cite{AW16}. 

One can now directly perturb the Fredholm section from the M-polyfold to the M-polyfold bundle as in~\cite{HWZ14} to make it transverse to zero. Alternatively, arguing as in~\cite{CMS03},
by finite-dimensional approximation one can write a compact zero-set of a Fredholm section 
from an M-polyfold to an M-polyfold bundle as 
the zero-set of a section from a finite-dimensional manifold to a finite-rank vector bundle over the manifold,
and it is well-known that such a section can be made transverse to zero by a small perturbation. 
The difference of the rank of the vector bundle and the dimension of the underlying manifold corresponds to the Fredholm index.
In the language of Cieliebak--Mundet--Salamon one can think of such a section  
as a finite-dimensional $G$-moduli problem for the trivial Lie group~$G$; 
in the language of Fukaya--Ono~\cite{FO99} and Fukaya--Oh--Ohta--Ono~\cite{FOOO09},
such a section corresponds to a global M-Kuranishi structure, 
namely a Kuranishi structure consisting of one chart with no orbifold flavour. 

When applying this approach to~$\RFH$, it is sufficient to work with one~$J$, that one is free to choose independent 
of~$t$ and~$\eta$.
In this approach the elements of the abstractly perturbed moduli spaces do not correspond to gradient flow lines 
anymore.
That the gradient flow lines with cascades of~$\RFH$ fit into the M-polyfolds set-up has not yet been worked out in detail.

\subsubsection*{A Conley-type argument.}
In the appendix we outline an intermediate approach to transversality, 
that is inspired by Conley index theory.
It allows for more general perturbations than the traditional approach,
namely non-local vector fields on the loop space~$\cL$, 
but still stays in the framework of gradient flow lines. 
Assume that one knows that $\mathcal{G}_a^b$ is compact for some~$J$.
The larger class of non-local perturbations then makes it easy to achieve transversality. 
The danger is now that even for arbitrarily small perturbations the spaces $\mathcal{G}_a^b(v)$
become non-compact. 
We shall show, however, that for sufficiently small perturbations~$v$ 
one can select a compact part of~$\mathcal{G}_a^b(v)$ near~$\mathcal{G}_a^b$:
For a flow line~$x$ of $-\nabla \cA +v$ let $\ev (x) = x(0) \in \cL$ be the evaluation at~$0$. 
Since $\mathcal{G}_a^b$ is compact, $K := \ev (\mathcal{G}_a^b) \subset \cL$ is compact. 
Choose a bounded open neighbourhood~$\cN$ of~$K$ in~$\cL$.
Then it turns out that for sufficiently small perturbations~$v$ the part $\mathcal{G}_a^b(v,\overline \cN)$ 
of~$\mathcal{G}_a^b(v)$ that lies in~$\overline \cN$ is compact and completely contained in~$\cN$,
that is, $\overline \cN$ is an isolating neighbourhood of the ``flow'' of~$-\nabla \cA +v$ for every small enough~$v$.
One can thus use the isolated invariant sets~$\mathcal{G}_a^b(v,\overline \cN)$  with regular~$v$ to define Floer homology. 
In Appendix~\ref{s:app.Morse} we describe this approach to transversality in detail 
in the framework of Morse homology on a finite-dimensional but possibly non-compact manifold.
The arguments are chosen in such a way that they translate to Floer homology (see~\S~\ref{ss:conley}), 
with the key difference that now the compactness of the selected components~$\mathcal{G}_a^b(v,\overline{\cN})$ 
does not just follow from the Arzel\`a--Ascoli theorem, 
because $\cL$ is not locally compact. 
To remedy for this, we prove in \S~\ref{ss:nonlocalcomp} a compactness result for non-local perturbations 
of the Cauchy--Riemann operator.

\m
\ni
{\bf Acknowledgments.}
We thank Alberto Abbondandolo for pointing out the reference~\cite{QuSa72}.
We are grateful to the referees for valuable comments and suggestions.

\section{Recollections on Rabinowitz--Floer homology} \label{s2}

In this section we recall the construction of the (non-equivariant) Rabinowitz--Floer homology 
of a hypersurface~$\Sigma$ of restricted contact type, following~\cite{CieFra09} and~\cite{AlbFra10}.
Our construction of equivariant Rabinowitz--Floer homology in the next section will be based on this construction. 

Consider an exact convex symplectic manifold $(M,\lambda)$.
This means that $\lambda$ is a one-form on the connected manifold~$M$ such that $d\lambda$ is a symplectic form,
and that $(M,d\lambda)$ is convex at infinity, i.e., 
there exists an exhaustion $M = \bigcup_k M_k$ of~$M$ by compact subsets $M_k \subset M_{k+1}$
with smooth boundaries~$\partial M_k$ such that $\lambda |_{\partial M_k}$ is a contact form.
We further fix a closed connected smooth hypersurface~$\Sigma$ in~$M$
that is bounding and of contact type.
The former means that $M \setminus \Sigma$ has two components, one compact and one non-compact, 
and the latter means that $\lambda |_\Sigma$ is a contact form, or equivalently 
that the vector field $Y_\lambda$ implicitly defined by $\iota_{Y_\lambda} d\lambda = \lambda$
is transverse to~$\Sigma$.

For a smooth function~$F$ on~$M$, the Hamiltonian vector field~$X_F$ is defined by $\iota_{X_F} d\lambda = dF$,
and $\varphi_F^t$ denotes the flow of~$X_F$.
The Reeb flow $\varphi_R^t$ on~$\Sigma$ is the flow of the vector field~$R$ defined by $d\lambda (R, \cdot) =0$ and $\lambda (R)=1$.

\subsection{The Rabinowitz action functional}

A {\it defining Hamiltonian for~$\Sigma$} is 
a smooth function $F \colon M \to \RR$ such that
$\Sigma \,=\, F^{-1}(0)$, such that $dF$ has compact support, 
and such that $\varphi_F^t$ restricts on~$\Sigma$ to the Reeb flow $\varphi_R^t$ of~$(\Sigma, \lambda |_\Sigma)$.
The set of defining Hamiltonians is non-empty and convex.
Given a defining Hamiltonian~$F$,
the Rabinowitz action functional
$\cA^F \colon \cL \times \RR \to \RR$ is defined by
\begin{equation}  \label{e:unperturbed}
\cA^F (v,\eta) \,=\, -\int_{S^1} v^*\lambda - \eta \int_{S^1} F \bigl(v(t) \bigr) \,dt .
\end{equation}
Its critical points $(v,\eta)$ are the solutions of the problem
$$
\dot v(t) = \eta \,X_F(v(t)), \quad 0 =  \int_{S^1} F(v(t)) \, dt ,
$$
i.e., pairs $(v,\eta)$
with $\eta \in \RR$ and $v$ a closed curve on~$\Sigma$ of the form $v(t) = \varphi_F^{\eta t}$, $t \in \RR$.
The critical points therefore correspond to closed orbits of~$X_F$ on the fixed energy surface 
$\Sigma = F^{-1}(0)$ of {\it arbitrary}\/ period $|\eta| \geqslant 0$.
\footnote{
Despite J.\ Moser's explicit statement that the action functional~\eqref{e:unperturbed} is useless
for finding periodic orbits, \cite[p.\ 731]{Mos76}, P.\ Rabinowitz in~\cite[p.\ 161 and (2.7)]{Rab78}
used precisely this functional to prove his celebrated existence theorem for periodic orbits 
on starshaped hypersurfaces in~$\RR^{2n}$,
thus pioneering the use of global critical point methods in Hamiltonian mechanics.
In \cite{CieFra09} and subsequent papers, the functional~\eqref{e:unperturbed} was therefore called 
Rabinowitz action functional.
Other good names for this functional may be 
``fixed energy action functional'' or ``Hamiltonian free period action functional'', 
since it selects solutions on the prescribed energy level $\{H=0\}$, 
allowing for arbitrary period~$|\eta|$.}
Since $v \subset \Sigma$ and $\varphi_F^t = \varphi_R^t$ along~$\Sigma$,
$$
\cA^F (v, \eta) \,=\, - \int_{S^1} v^* \lambda \,=\, -\eta ,
$$
that is, the critical values of $\cA^F$ are zero and minus the periods of the closed Reeb orbits on~$\Sigma$.

The action functional $\cA^F$ is invariant under the $S^1$-action on $\cL \times \RR$ given by
\begin{equation} \label{e:s1.action}
\tau \cdot \bigl( v(\cdot), \eta \bigr) \,\mapsto\, \bigl( v(\cdot -\tau), \eta \bigr) .
\end{equation}
Therefore, the functional $\cA^F$ is never Morse.
The component $\{ (p,0) \mid p \in \Sigma\} \cong \Sigma$ of the critical set is always Morse--Bott for $\cA^F$,
see \cite[Lemma~2.12]{AlbFra10}.
The following assumption on~$\Sigma$ is sufficient for $\cA^F$ to be Morse--Bott:
\begin{equation} \label{e:ass}
\mbox{\it Every periodic orbit of the Reeb flow $\varphi_R^t$ is non-degenerate.}
\end{equation}
In other words, for a $T$-periodic orbit $\gamma$ of the Reeb flow, 
$1$ is not in the spectrum of the linearization $T_p\varphi_R^T \colon \xi_p \to \xi_p$ at $p=\gamma (0)$,
where $\xi = \ker \lambda$ denotes the contact structure of~$\Sigma$. 
This holds if and only if for any defining Hamiltonian~$F$ of~$\Sigma$,
for every periodic orbit of $\varphi_F^t$ on~$\Sigma$ the Floquet multiplier~$1$ has multiplicity~$2$.

\subsection{Rabinowitz--Floer homology}  \label{ss:RFH}

Rabinowitz--Floer homology $\RFH (\Sigma, M)$ is the Floer homology of the functional~$\cA^F$, 
where $F$ is any defining Hamiltonian for~$\Sigma$.
We assume the reader to be familiar with the construction in~\cite{CieFra09},
and also refer to~\cite{AlbFra10} and to the survey~\cite{AlbFra12}.
Here, we only point out a few aspects in the construction of~$\RFH (\Sigma, M)$ 
that do not arise in the construction of usual Hamiltonian Floer homology.

\medskip
{\bf 1. The chain groups.}
The functional $\cA^F$ is not Morse, but Morse--Bott. 
One therefore chooses an auxiliary Morse function $h \colon \Crit \cA^F \to \RR$,
and generates the chain groups by the critical points of~$h$.
However, even though the symplectic form $d \lambda$ is exact, the generators of the Rabinowitz--Floer chain groups
$\FC(\cA^F, h)$ are not finite sums $\sum \xi_c c$ with $\xi_c \in \ZZ_2$ and $c \in \Crit h$, 
but possibly infinite sums $\sum \xi_c c$ that for every $\kappa \in \RR$ satisfy the finiteness condition
$$
\# \left\{ c \in \Crit h \mid \xi_c \neq 0, \, \cA^F(c) \geqslant \kappa \right\} \,<\, \infty .
$$
This must be done so for the following reason:
Assume that $c$ lies on the critical point $(v, \eta)$ of~$\cA^F$, with $\eta \neq 0$.
Then $\cA^F (v, \eta) = -\eta$.
Since with $(v, \eta)$ also $(v, k\eta)$ belongs to $\Crit \cA^F$ for each $k \in \ZZ$,
we see that $\cA^F$ is not bounded from below on $\Crit \cA^F$.
Hence there may be infinitely many critical points that appear in the image $\pp c$ of the boundary operator.

\smallskip
{\bf 2. The almost complex structures.}
Let $\cJ_{\con}$ be the set of almost complex structures on~$M$ that are $d\lambda$-compatible and 
convex at infinity.
The choice of the set of almost complex structures used to define $\RFH (\Sigma, M)$ depends on the
method that one uses to establish transversality (cf.\ the introduction).
If one works with polyfolds or with the Conley-type approach explained in the appendix, 
one can take a fixed $J \in \cJ_{\con}$.
In the next paragraph we describe the boundary operator in the traditional way.
For this we fix $J_* \in \cJ_{\con}$ and following~\cite{AbbMer15} consider the set $\mathcal{J}$ of 
smooth $S^1 \times \RR$-families $\bJ = \{J_t (\cdot, \eta)\} \subset \cJ_{\con}$
such that
\begin{equation} \label{e:boundJc}
\sup_{t,\eta} \| J_t (\cdot, \eta) \|_{C^\ell} \,<\, \infty \quad \mbox{ for all }\, \ell \in \NN
\end{equation}
and such that there exists a constant $c>1$ (depending on the family) for which
\begin{equation} \label{e:boundAM}
\tfrac 1c \, \| J_*(x) \| \,\leqslant\, \| J_t(x,\eta) \| \,\leqslant\, c \, \| J_*(x) \|
\quad \mbox{ for all }\, x \in M \mbox{ and } (t,\eta) \in S^1 \times \RR .
\end{equation}
Here, $\| \cdot \|$ is the norm taken with respect to some background Riemannian metric on~$M$.

\smallskip
{\bf 3. The boundary operator.}
The boundary operator $\partial$ on $\FC (\cA^F, h)$ is defined by counting gradient flow lines with cascades
(see \cite[Appendix~A]{Fra04}).
These flow lines consist of (partial) negative gradient flow lines of~$h$
and finite energy Floer gradient flow lines of~$\cA^F$.
Given a family $\bJ \in \mathcal{J}$ and two critical points $(v_-,\eta_-)$ and $(v_+,\eta_+)$ of~$\cA^F$, 
a Floer gradient flow line is a solution $(v,\eta) \in C^\infty (\RR \times S^1, M \times \RR)$ of the problem
\begin{equation} \label{e:Floertra}
\left. 
\begin{array}{rcl}
\pp_s v (s,t) + J_t\bigl(v(s,t),\eta(s)\bigr) \bigl( \pp_t v(s,t) - \eta (s) \,X_F \bigl( v(t) \bigr) &=& 0 , \\ [0.4em]
\dot \eta (s) + \int_{S^1} F \bigl(v(s,t) \bigr) \,dt &=& 0 ,
\end{array}
\right\}
\end{equation}
with asymptotic boundary conditions~$(v_-,\eta_-)$ and $(v_+,\eta_+)$.
The main analytical issue in defining the boundary operator~$\partial$ 
is to prove a uniform $L^\infty$-bound on the $\eta$-component of the solutions of~\eqref{e:Floertra} 
with given boundary conditions.
This is done in~\cite[Corollary~3.3]{CieFra09} for $\eta$-independent~$J$, 
and the proof goes through thanks to~\eqref{e:boundAM}.
Assumption~\eqref{e:boundJc} is imposed to avoid bubbling, so that the space of all solutions
of~\eqref{e:Floertra} is $C^\infty_{\loc}$-compact.
Transversality for the space of solutions of~\eqref{e:Floertra} between two critical points
for a generic set of~$\bJ \in \mathcal{J}$ is proven in~\cite[\S~4.3]{AbbMer15}.

We remark that
the construction of the boundary operator by gradient flow lines with cascades in~\cite[Appendix~A]{Fra04} 
is given for Morse homology on finite-dimensional manifolds. 
While this construction directly carries over to the case of Floer homology, 
some parts of this generalisation (such as gluing) are not worked out in the literature.
The same applies to the $S^1$-equivariant Rabinowitz--Floer homology described in the next section.
The foundational work coming closest to the holomorphic curve set-up considered in this paper 
is in~\cite{Bou02, BouOan09:Duke} and~\cite[\S 10]{Schm16}.
Another way to rigorously establish~$\RFH$ and~$\RFH^{S^1}$ is by verifying that the flow lines with cascades
fit into the $M$-polyfold set-up (cf.\ the introduction).

\smallskip
{\bf 4. Invariance.}
The resulting homology group $\FH (\cA^F) := \ker \pp / \im \pp$ does not depend on the choice of a defining 
function~$F$ for~$\Sigma$.
One can therefore define $\RFH (\Sigma, M) := \FH (\cA^F)$ for any choice of~$F$.
Moreover, given two bounding contact hypersurfaces $\Sigma_0$ and~$\Sigma_1$ that are isotopic
through a family~$\{ \Sigma_s \}_{0 \leqslant s \leqslant 1}$ of contact hypersurfaces, 
\begin{equation} \label{e:iso01}
\RFH (\Sigma_0, M) \,\cong\, \RFH (\Sigma_1, M) .
\end{equation}
For the proof, one chooses a smooth family $F_s \colon M \to \RR$ of defining Hamiltonians for~$\Sigma_s$
such that $F_s = F_0$ for $s \leqslant 0$ and $F_s=F_1$ for $s \geqslant 1$,
and uses solutions of~\eqref{e:Floertra} with $F$ replaced by~$F_s$ to construct 
a chain homotopy equivalence 
between $\FC (\cA^{F_0}, h_0)$ and $\FC (\cA^{F_1}, h_1)$.
The main analytical issue is again proving a bound on the $\eta$-components, 
which can be done as in \cite[Corollary~3.4]{CieFra09} thanks to~\eqref{e:boundAM}.

Recall that we have worked for now under the assumption~\eqref{e:ass}.
This assumption on~$\Sigma$ is generic in the $C^\infty$-topology.
In view of~\eqref{e:iso01} we can define the Rabinowitz--Floer homology $\RFH (\Sigma, M)$
of any bounding contact hypersurface as $\RFH (\Sigma', M)$ where $\Sigma'$ is a close-by hypersurface
meeting assumption~\eqref{e:ass}.


\section{Construction of equivariant Rabinowitz--Floer homology} \label{s3}

In this section we give a Borel-type construction of $S^1$-equivariant Rabinowitz--Floer homology,
closely following the
construction of $S^1$-equivariant symplectic homology given by Viterbo in~\cite[\S 5]{Vit99},
see also~\cite{BouOan09:Gysin}.

\subsection{The equivariant Rabinowitz action functional}

For each integer $N \geqslant 1$ denote by $S^{2N+1}$ the odd-dimensional unit sphere in~$\CC^{N+1}$. 
The circle $S^1$ acts on $S^{2N+1}$ by
$$
\tau \cdot (z_1, \dots, z_{N+1}) \,=\, (\tau z_1, \dots, \tau z_{N+1}) .
$$ 
The quotient of this action is complex projective space $\CC P^N = S^{2N+1}/S^1$.
Recall the action~\eqref{e:s1.action} of $S^1$ on the loop space $\cL$, and let $S^1$ act on
$\cL \times \RR \times S^{2N+1}$ by the diagonal action
\begin{equation} \label{def:diag.action}
\tau \cdot \bigl( v(\cdot), \eta, z \bigr) \,=\, \bigl( v(\cdot - \tau), \eta, \tau \cdot z \bigr) .
\end{equation}
We shall denote the circle~$S^1$ with this action on $\cL \times \RR \times S^{2N+1}$ by~$\TT$.
Denote the quotient of this action by $\cL \times \RR \times_{\TT} S^{2N+1}$.
The functional 
$\widetilde \cA^{F,N;\TT} \colon \cL \times \RR \times S^{2N+1} \to \RR$
defined by 
\begin{equation}  \label{e:unperturbed.times}
\widetilde \cA^{F,N;\TT} (v,\eta,z) \,=\, -\int_{S^1} v^*\lambda - \eta \int_{S^1} F \bigl(v(t) \bigr) \,dt 
\end{equation}
is Morse--Bott if and only if the functional~$\cA^{F}$ defined in~\eqref{e:unperturbed} is Morse--Bott.
Indeed, the critical set of $\widetilde \cA^{F,N;\TT}$ is the critical set of~$\cA^{F}$ times $S^{2N+1}$.
Since the functional~\eqref{e:unperturbed.times} is invariant under the action~\eqref{def:diag.action},
we can define the {\it equivariant Rabinowitz action functional}\/
$\cA^{F,N;\TT} \colon \cL \times \RR \times_\TT S^{2N+1} \to \RR$
by 
\begin{equation}  \label{e:unperturbed.S1}
\cA^{F,N;\TT} ([v,\eta,z]) \,=\, -\int_{S^1} v^*\lambda - \eta \int_{S^1} F \bigl(v(t) \bigr) \,dt ,
\end{equation}
and since the action~\eqref{def:diag.action} is free, 
this functional is Morse--Bott under the assumption~\eqref{e:ass} on~$\Sigma$.

\subsection{Equivariant Rabinowitz--Floer homology}  \label{ss:eqRFH}

$\TT$-equivariant Rabinowitz--Floer homology $\RFH^\TT (\Sigma, M)$ is the direct limit in~$N$ 
of the Floer homology of the functional~$\cA^{F,N;\TT}$, 
where $F$ is any defining Hamiltonian for~$\Sigma$.

\medskip
{\bf 1. The chain groups.}
Fix a defining Hamiltonian~$F$ for~$\Sigma$ meeting assumption~\eqref{e:ass},
and fix $N \in \NN$. Then $\widetilde \cA^{F,N;\TT}$ is Morse--Bott, with critical manifolds
the union of $\Sigma \times \{0\} \times S^{2N+1}$ and $C_i \times \{k \2 \eta_i\} \times S^{2N+1}$,
$k \in \ZZ \setminus \{0\}$,
where each $C_i \times \{\eta_i\}$ is a circle of simple Reeb orbits of period~$\eta_i$.
Since the action of~$\TT$ on $\cL \times \RR \times S^{2N+1}$ is free, 
$$
\Crit \cA^{F,N;\TT} \,=\, \Crit \widetilde \cA^{F,N;\TT} / \TT \,=\, \Crit \cA^F \times_\TT S^{2N+1}
$$
is a closed manifold.
Denote by~$g_{S^{2N+1}}$ the round Riemannian metric on~$S^{2N+1}$, and choose a Riemannian metric~$g_\Sigma$ on~$\Sigma$
and $S^1$-invariant Riemannian metrics $g_{C_i}$ on~$C_i$.
Then the Riemannian metric~$g_N$ on~$\Crit \widetilde \cA^{F,N;\TT}$ defined by
$g_N |_{\Sigma \times \{0\} \times S^{2N+1}} = g_\Sigma \oplus g_{S^{2N+1}}$ and 
$g_N |_{C_i \times \{ k \2 \eta_i \} \times S^{2N+1}} = g_{C_i} \oplus g_{S^{2N+1}}$
is $\TT$-invariant,
and hence descends to the Riemannian metric $g_N^\TT$ on~$\Crit \cA^{F,N;\TT}$.
Choose a Morse function $h_N \colon \Crit \cA^{F,N;\TT} \to \RR$ such that the pair $(h_N,g_N^\TT)$
is Morse--Smale 
(that is, the stable and unstable manifolds of the negative gradient flow of~$h_N$ with respect to~$g_N^\TT$ intersect transversally).
The chain group $\FC (\cA^{F,N;\TT},h_N)$ consists of Novikov sums $\sum \xi_c c$ with $c \in \Crit h_N$, as in Section~\ref{ss:RFH}.

\smallskip
{\bf 2. The almost complex structures.}
If one works with polyfolds or with the Conley-type approach, 
one can, again, just take a fixed $J \in \cJ_{\con}$.
Here, we again fix $J_* \in \cJ_{\con}$ and look at 
smooth $S^1 \times S^{2N+1} \times \RR$-families $\bJ = \{ J_{t,z} (\cdot, \eta)\} \subset \cJ_{\con}$
such that
\begin{equation} \label{e:boundJcz}
\sup_{t,z,\eta} \| J_{t,z} (\cdot, \eta) \|_{C^\ell} \,<\, \infty 
\quad \mbox{ for all }\, \ell \in \NN 
\end{equation}
and such that there exists a constant $c>1$ (depending on the family) for which
\begin{equation} \label{e:boundAMz}
\tfrac 1c \, \| J_*(x) \| \,\leqslant\, \| J_{t,z}(x,\eta) \| \,\leqslant\, c \, \| J_*(x) \|
\quad \mbox{ for all }\, x \in M \mbox{ and } (t,z,\eta) \in S^1 \times S^{2N+1} \times \RR .
\end{equation}
Furthermore, we impose that the family~$\bJ$ is $S^1$-invariant:
\begin{equation} \label{e:Jinv}
J_{t+\tau, \tau z} (\cdot, \eta) \,=\, J_{t,z} (\cdot, \eta) \quad\, 
\mbox{ for all $(t,z,\eta) \in S^1 \times S^{2N+1} \times \RR$ and $\tau \in S^1$.}
\end{equation}
The space $\mathcal{J}^{S^1}$ of all families $\bJ$ in~$\cJ_\con$
satisfying~\eqref{e:boundJcz}, \eqref{e:boundAMz} and~\eqref{e:Jinv}
is non-empty (since property~\eqref{e:Jinv} is obtained by averaging over~$S^1$) 
and contractible.

\smallskip
{\bf 3. The boundary operator.}
Let $\widetilde h_N \colon \Crit \widetilde \cA^{F,N;\TT} \to \RR$
be the lift of~$h_N$.
Then $\widetilde h_N$ is Morse--Bott, with $\TT$-orbits as critical manifolds.
Given two critical points $c^+$, $c^-$ of~$h_N$,
denote by~$C^+$, $C^-$ the corresponding critical circles of~$\widetilde h_N$.
Given $\bJ \in \mathcal{J}^{S^1}$
consider all gradient flow lines with cascades $\widehat \cM (c^+,c^-)$ from a point in~$C^+$ to a point in~$C^-$.
Here, the (partial) Morse flow lines are (partial) negative gradient flow lines of~$\widetilde h_N$ 
on~$\Crit \widetilde \cA^{F,N;\TT}$ with respect to~$g_N$, and 
the cascades (i.e., the Floer gradient flow lines) are finite energy solutions 
$(v,\eta,z) \in C^\infty(\RR \times S^1, M \times \RR \times S^{2N+1})$
of the problem
\begin{equation} \label{e:Floertra.S1}
\left. 
\begin{array}{rcl}
\pp_s v (s,t) + J_{t,z(s)}\bigl(v(s,t),\eta(s) \bigr) \bigl( \pp_t v(s,t) - \eta (s) \,X_F \bigl( v(t) \bigr) &=& 0 , \\ [0.4em]
\dot \eta (s) + \int_{S^1} F \bigl(v(s,t) \bigr) \,dt &=& 0 , \\ [0.4em]
\dot z (s) + \nabla_{g_{S^{2N+1}}} \widetilde h_N (z(s))       &=& 0 .
\end{array}
\right\}
\end{equation}
Here, $\nabla_{g_{S^{2N+1}}} \widetilde h_N (z)$ denotes the component of $\nabla_{g_N} \widetilde h_N (z)$ along $T_zS^{2N+1}$.
Since $g_N$ and~$J$ are $\TT$-invariant, $\TT$ freely acts on $\widehat \cM (c^+,c^-)$.
The space $\widehat \cM (c^+,c^-)$ therefore decomposes as 
$$
\widehat \cM (c^+,c^-) \,=\, \coprod_{c \in C^+} \widehat \cM (c,c^-) 
$$
where $\widehat \cM (c,c^-)$ is the space of gradient flow lines with cascades from $c \in C^+$ 
with the last gradient flow line of~$\widetilde h_N$ converging to an arbitrary point in~$C^-$,
and $\widehat \cM (c^+,c^-) / \TT \cong \widehat \cM (c,c^-)$ for any $c \in C^+$.
One shows as in~\cite[\S~4.3]{AbbMer15} that for a generic subset of families $\bJ \in \mathcal{J}^{S^1}$ 
the spaces~$\widehat \cM (c,c^-)$ are smooth manifolds. 

The real numbers $s \in \RR$ freely act by shift on each Floer gradient flow line in a gradient flow line with cascades 
in~$\widehat \cM (c^+,c^-)$.
The space $\cM (c^+,c^-) \cong \coprod_{c \in C^+} \cM (c,c^-)$ obtained by modding out these $\RR$-actions is compact.
The main point in the proof is, again, a uniform $L^\infty$-bound on the $\eta$-component of the solutions of~\eqref{e:Floertra.S1} 
with given boundary conditions. Such a bound is obtained exactly as in~\cite[Corollary~3.3]{CieFra09},
thanks to~\eqref{e:boundAMz}.

Now the boundary operator on $\FC (\cA^{F,N;\TT},h_N)$ is defined by 
$$
\pp \1 (c^+) \,=\, \sum_{c^-} \nu (c^+,c^-) \,c^-
$$
where the sum runs over those $c^-$ for which $\cM (c^+,c^-) / \TT \cong \cM (c,c^-)$ is $0$-dimensional and where
$\nu (c^+,c^-)$ is the number mod~$2$ of elements in this space.

\smallskip
{\bf 4. Invariance.}
Let $\FH (\cA^{F,N;\TT},h_N,J) := \ker \pp / \im \pp$ be the resulting homology groups.
The inclusion $S^{2N+1} \to S^{2N+3}$ is $\TT$-equivariant.
In particular, $\Crit \cA^{F,N;\TT} \subset \Crit \cA^{F,N+1;\TT}$.
Since $g_{S^{2N+3}}$ restricts to~$g_{S^{2N+1}}$ on~$S^{2N+1}$, the Riemannian metric~$g_{N+1}$ restricts to~$g_N$ 
on~$\Crit \widetilde \cA^{F,N;\TT}$.
Given a Morse function $h_N$ on~$\Crit \cA^{F,N;\TT}$ as above, we choose a Morse function $h_{N+1}$ on~$\Crit \cA^{F,N+1;\TT}$
such that $h_{N+1}$ extends $h_N$, such that $\Crit h_N \subset \Crit h_{N+1}$, 
and such that the pair $(h_{N+1}, g_{N+1}^\TT)$ is Morse--Smale.
Further, we choose the family $J_{N+1} = J_{t,z}(\cdot, \eta)$ with $z \in S^{2N+3}$ such that it extends the family 
$J_N = J_{t,z}(\cdot, \eta)$ with $z \in S^{2N+1}$.
The chain complex $\FC (\cA^{F,N;\TT},h_N,J_N)$ is thus a subcomplex of $\FC (\cA^{F,N+1;\TT},h_{N+1},J_{N+1})$.
We thus obtain a homomorphism 
\begin{equation} \label{e:dir}
\iota_N \colon \FH (\cA^{F,N;\TT},h_N,J_N) \to \FH (\cA^{F,N+1;\TT},h_{N+1},J_{N+1}) .
\end{equation}
The groups $\FH (\cA^{F,N;\TT},h_N,J_N)$ do not depend on the choice of~$h_N$ and~$J_N$,
nor on the choice of~$g_\Sigma$ in the definition of~$g_N$,
nor on the defining Hamiltonian~$F$ for~$\Sigma$. 
This is proven by Floer continuation as in~\cite{CieFra09} (see also Section~\ref{s:invariance}).
These continuation isomorphisms commute with the inclusion homomorphisms in~\eqref{e:dir}:
Given another defining Hamiltonian~$F'$ and other choices~$h_N'$ and~$J_N'$,
there is a commutative diagram
$$
\xymatrixcolsep{4pc}\xymatrix{ 
\FH (\cA^{F,N;\TT},h_N,J_N)
\ar[d]^{\cong} \ar[r]^-{\iota_N} &
\FH (\cA^{F,N+1;\TT},h_{N+1},J_{N+1}) \ar[d]^{\cong} & 
\\
\FH (\cA^{F',N;\TT},h_N',J_N') \ar[r]^-{\iota_N'} & 
\FH (\cA^{F',N+1;\TT},h_{N+1}',J_{N+1}')  .
}
$$
The direct limit
\begin{equation}  \label{e:dirlim}
\RFH^\TT (\Sigma, M) := \varinjlim \FH (\cA^{F,N;\TT},h_N,J_N) 
\end{equation}
therefore only depends on $\Sigma$. 
In fact,  $\RFH^\TT (\Sigma, M)$ is invariant under isotopies of 
bounding contact hypersurfaces (cf.\ Section~\ref{ss:RFH}).

\begin{remarks} \label{rem:grading}
{\rm
{\bf 1.}
Our homology groups $\RFH^\TT (\Sigma, M)$ are not graded.
We therefore do not need to assume that the first Chern class of~$(M,d\lambda)$ vanishes on~$\pi_2(M)$.
Under this assumption, the groups $\RFH^\TT (\Sigma, M)$ carry a $\ZZ$-grading (with values in $\frac 12 + \ZZ$),
cf.\ \cite[Section~4]{CieFra09}.

\smallskip
{\bf 2.}
The above construction of $S^1$-equivariant Rabinowitz--Floer homology should give the same result 
as the construction in~\cite{BouOan09:Gysin} 
which uses parametrized symplectic homology,
when applied to the parameter space~$\RR \times S^{2N+1}$: 
The difference in the construction is that our parameter space $\RR \times S^{2N+1}$ is not compact,
and that we work with cascades instead of suitable perturbations of the Hamiltonian~$F$.  
We expect that 
combining the construction in~\cite{BouOan09:Gysin} with the $L^\infty$-estimates on the $\eta$-component 
from~\cite[Section~3]{CieFra09}
leads to the same groups~$\RFH^\TT (\Sigma, M)$ 
in view of a version of the Correspondence Theorem~3.7 in~\cite{BouOan09:Duke}.

A construction of an $S^1$-equivariant Rabinowitz--Floer homology that stays within the setting of 
$S^1$-equivariant symplectic homology was given recently in~\cite{CieOan15}.
We expect that also this homology is isomorphic to~$\RFH^\TT (\Sigma, M)$.
}
\end{remarks}

\section{Proof of Theorem A}  \label{s:proofThA}

In this section we prove our main result: $\RFH^\TT(\Sigma, M)=0$ if $\Sigma$ is displaceable.
For the proof, we first recall how the analogous result is proven in the non-equivariant case.
We shall apply the same method in the non-equivariant case.

\subsection{The perturbed Rabinowitz action functional, and leafwise intersections}

It has been shown in~\cite{CieFra09} that $\RFH (\Sigma, M)$ 
vanishes if $\Sigma$ is displaceable.
This result has been reproved in~\cite{AlbFra10} by a more geometric argument,
in which the functional~$\cA^F$ is perturbed to a functional whose critical points
are leafwise intersections.
While the argument in~\cite{CieFra09} can be useful in problems where 
the leafwise intersection argument does not help
(such as proving the existence of a closed characteristic on a displaceable stable hypersurface~\cite{CieFraPat10}),
we here apply the leafwise intersection argument from~\cite{AlbFra10}.


A {\it perturbation pair}\/ for the Rabinowitz action functional is a tuple
$$
(\chi,H)  \in  C^\infty \bigl(S^1,[0,\infty) \bigr) \times C^\infty \bigl( M \times S^1, \RR \bigr)
$$ 
such that
$\int_{S^1} \chi(t) \,dt \,=\,1$.
For a perturbation pair, the {\it perturbed Rabinowitz action functional}\/
$\cA^F_{\chi,H} \colon \cL \times \RR \to \RR$ 
is defined by
\begin{equation} \label{e:per.free}
\cA^F_{\chi,H} (v,\eta) \,=\, -\int_{S^1} v^*\lambda - \eta \int_{S^1}
\chi(t) \, F\bigl(v(t)\bigr) \,dt - \int_{S^1} H\bigl(v(t),t \bigr) \,dt.
\end{equation}
The critical points $(v,\eta)$ of this perturbed action functional are the solutions of 
the system
\begin{equation} \label{li}
\left. 
\begin{array}{rcl}
\dot v(t) &=& \eta \, \chi(t) \,X_F \bigl( v(t) \bigr) + X_H \bigl( v(t),t \bigr), \\ [0.4em]
0 &=& \int_{S^1} \chi(t) \,F \bigl(v(t) \bigr) \,dt .
\end{array}
\right\}
\end{equation}

As noticed in \cite{AlbFra10}, it is useful to look at special perturbation pairs:

\begin{definition} \label{def:Moser}
{\rm 
A perturbation pair $(\chi,H)$ is called of {\it Moser type}\/ if
there exists $t_0 \in S^1$ such that the time support of~$H$ lies in
$[t_0,t_0+1/2]$ and the support of $\chi$ lies in $[t_0-1/2,t_0]$.
}
\end{definition}

The energy hypersurface $\Sigma = F^{-1}(0)$ is foliated by its leaves
$L_x = \left\{ \varphi_F^t(x) \mid t \in \RR \right\}$.
Given a perturbation $H$ as above, a point $x \in \Sigma$ is called a 
{\it leafwise intersection point for~$H$}
if $\varphi_H^1(x) \in L_x$.
The following lemma was observed in~\cite{AlbFra10}.

\begin{lemma} \label{l:intersection}
If a perturbation pair is of Moser type and $(v,\eta)$ is a solution of~\eqref{li}, 
then $v(t_0)$ is a leafwise intersection point for~$H$ on~$\Sigma = F^{-1}(0)$. 
\end{lemma}


\subsection{The perturbed equivariant Rabinowitz action functional}

In order to show that $\RFH^\TT (\Sigma, M)$ vanishes for displaceable $\Sigma$,
we wish to apply the same method as in the non-equivariant case.

In the following $S^1$ acts diagonally on $S^1 \times S^{2N+1}$ by $\tau( \cdot , z) = ( \cdot - \tau, \tau z)$,
and $S^1 \times_{S^1} S^{2N+1}$ is the quotient of $S^1 \times
S^{2N+1}$ under this action. 
A {\it perturbation triple}\/ is a triple
$$
(\psi,G,k) \in 
C^\infty \bigl( S^1 \times_{S^1} S^{2N+1},[0,\infty) \bigr) 
\times
C^\infty \bigl( M \times S^1 \times_{S^1} S^{2N+1}, \RR \bigr)
\times
C^\infty \bigl( \CC P^N, \RR \bigr)
$$
such that for every $z \in S^{2N+1}$,
\begin{equation} \label{e:intpsi1}
\int_{S^1} \psi \bigl( [t,z] \bigr) \,dt \,=\,1 ,
\end{equation}
and such that $k$ is a Morse function on~$\CC P^N$.
For a perturbation triple we define the perturbed equivariant Rabinowitz action functional
\begin{equation} \label{def:cA.per}
\cA_{\psi,G,k} \,:=\, \cA^{F,N;\TT}_{\psi,G,k}
\colon \cL \times \RR \times_{\TT} S^{2N+1} \to \RR
\end{equation}
by
$$
\cA_{\psi,G,k} \bigl( [v,\eta,z] \bigr) \,=\, 
-\int_{S^1} v^*\lambda 
-\eta \int_{S^1} \psi \bigl( [t,z] \bigr) F \bigl( v(t) \bigr) \,dt 
-\int_{S^1} G \bigl( v(t),[t,z] \bigr) \, dt - k \bigl( [z] \bigr) .
$$ 
Denote by
\begin{eqnarray*}
\begin{array}{lcl}
\widetilde \psi &\in&  C^\infty \bigl( S^1 \times S^{2N+1}, [0,\infty) \bigr), \\ [0.2em]
\widetilde G    &\in&  C^\infty (M \times S^1 \times S^{2N+1},\RR), \\ [0.2em]
\widetilde k    &\in&  C^\infty (S^{2N+1},\RR) 
\end{array}
\end{eqnarray*}
the lifts of $\psi$, $G$ and~$k$.
We can then write the lift of $\cA_{\psi,G,k}$ to $\cL \times \RR \times S^{2N+1}$ as
\begin{equation} \label{def:cA.per.lift}
\cA_{\widetilde \psi, \widetilde G, \widetilde k} \left(v,\eta,z \right) \,=\, 
-\int_{S^1} v^*\lambda 
-\eta \int_{S^1} \widetilde \psi \left(t,z \right) F \bigl( v(t) \bigr) \,dt 
-\int_{S^1} \widetilde G \bigl( v(t),t,z \bigr) \, dt - \widetilde k(z) .
\end{equation}
The critical points $(v,\eta,z)$ of~$\cA_{\widetilde \psi, \widetilde G, \widetilde k}$ are the solutions of the system
\begin{equation} \label{crit.tilde}
\left.
\begin{array}{rcl}
\dot v(t) &=& \eta \, \widetilde\psi(t,z) \,X_F \bigl( v(t) \bigr) + X_{\widetilde G} \bigl( v(t),t \bigr), \\ [.4em]
0 &=& \int_{S^1} \widetilde \psi(t,z)\, F \bigl(v(t) \bigr) \,dt , \\ [.4em]
0 &=& \eta \int_{S^1} F \bigl( v(t) \bigr) \2 \pp_z \widetilde \psi (t,z)\,dt - \int_{S^1} \pp_z \widetilde G \bigl( v(t), t,z)\, dt - d\widetilde k(z) .
\end{array}
\right\}
\end{equation}

\begin{definition} \label{def.adm}
{\rm 
A perturbation triple $(\psi,G,k)$ is called \emph{admissible} if
the following two conditions hold.

\begin{itemize}
\item[(i)] 
For each $z \in S^{2N+1}$ and each solution $(v,\eta)$ of equation~\eqref{li} 
with respect to the perturbation $(\widetilde \psi_z,\widetilde G_z)$
the identity $F \bigl( v(t) \bigr) \,d\widetilde \psi_t(z)=0$ 
holds for all $t \in S^1$.

\m
\item[(ii)] 
$|d \widetilde G_{x,t}(z) \,\hat z |< | d \widetilde k (z) \, \hat{z}|$
\quad for all\; $z \notin \Crit \widetilde k$, \;
$\hat z \neq 0 \in T_z S^{2N+1}$, \; $(x,t) \in M \times S^1$.
\end{itemize}
}
\end{definition}

\begin{lemma}  \label{l:crit}
Assume that $(\psi,G,k)$ is an admissible perturbation triple. 
Then critical points $[v,\eta,z]$ of $\cA_{\psi,G,k}$ have the property that 
$[z]$ is a critical point of~$k$, and for each $z \in S^{2N+1}$ over $[z]$,
the pair $(v,\eta)$ is a solution to equation~\eqref{li} for the perturbation
$(\widetilde \psi_z,\widetilde G_z)$.
\end{lemma}

\proof 
In view of the first two equations in~\eqref{crit.tilde}, we see that $(v,\eta)$
is a solution of~\eqref{li} for the perturbation
$(\widetilde \psi_z,\widetilde G_z)$. 
It remains to show that $[z]$ is a critical point of~$k$. 
In view of the last equation in~\eqref{crit.tilde}, 
for every $\hat{z} \in T_z S^{2N+1}$ the equation
$$
\eta \int_{S^1} F \bigl( v(t) \bigr) \, d \widetilde \psi_t(z) \,\hat{z} \,dt + 
\int_{S^1} d \widetilde G_{v,t} (z) \, \hat{z} \,dt + d \widetilde k (z) \,\hat{z} \,=\, 0
$$ 
has to be met. 
By assertion~(i) of Definition~\ref{def.adm},
the first term vanishes. 
Now assertion~(ii) implies $d \widetilde k(z) \,\hat{z}=0$, 
hence $[z]$ is a critical point of~$k$. 
\proofend

\begin{definition} \label{def.equiv.extension}
{\rm
Given a perturbation pair of Moser type $(\chi,H)$, we call a
perturbation triple $(\psi,G,k)$ an \emph{equivariant extension} of 
$(\chi,H)$ if the following conditions hold.
\begin{itemize}
\item[(I)] 
The perturbation triple $(\psi,G,k)$ is admissible. 
 
\m
\item[(II)] 
For every $z \in \Crit \widetilde k$ there exists $t_z \in S^1$ such that 
for every $t \in S^1$ and every $x \in M$ the identities $\widetilde G(x,t,z)=H(x,t+t_z)$ 
and $\widetilde \psi (t,z) = \chi(t+t_z)$ hold true.
\end{itemize}
}
\end{definition}


\begin{lemma} \label{l:extension}
For any perturbation pair $(\chi,H)$ of Moser type, 
there exists an equivariant extension.
\end{lemma}

\proof
Choose a Morse function $k$ on $\CC P^N$. 
For every $y \in \Crit k$ choose open neighborhoods
$$
y \in U_y \subset \overline{U_y} \subset V_y \subset \overline{V_y} \subset W_y
$$ 
with the property that $W_y$
is contractible, and for different critical points $y$ and~$y'$ of~$k$ 
the neighborhoods $W_y$ and $W_{y'}$ are disjoint. 
Since $W_y$ is contractible, the principal $S^1$-bundle 
$\pi \colon S^{2N+1} \to \CC P^N$ can be trivialized over~$W_y$. 
We abbreviate
$$
X \,=\, \bigcup_{y \in \Crit k} \pi^{-1}(W_y)
$$
and choose a trivialization
$$\Phi \colon X \to \pi(X) \times S^1.$$
We further choose smooth cutoff functions $\beta_1, \beta_2 \in
C^\infty (\CC P^N,[0,1])$ with the property that for every $y \in \Crit k$,
$$
\beta_1 |_{U_y}=1, \quad \beta_2|_{V_y}=1
$$
and
$$
\supp (\beta_1) \subset \bigcup_{y \in \Crit k} V_y, \qquad
\supp (\beta_2) \subset \bigcup_{y \in \Crit k} W_y .
$$ 
We further abbreviate by
$p \colon X \times S^1 \to S^1$
the projection to the second factor. We now set
$$
\widetilde G (x,t,z) \,=\, 
\left\{
\begin{array}{ll}
\beta_1([z]) \, H \bigl( x,t+p(\Phi(z)) \bigr), & z \in X, 
\vspace{0.2em}\\
0, & z \notin X.
\end{array}
\right.
$$
and
$$
\widetilde \psi (t,z) \,=\, 
\left\{
\begin{array}{ll}
\beta_2([z]) \, \chi \bigl( t+p (\Phi(z)) \bigr) + 1 - \beta_2([z]), & z \in X, 
\vspace{0.2em}\\ 
1, & z \notin X. 
\end{array}
\right.
$$ 
Define $G$ and $\psi$ by $G(x,[t,z]) = \widetilde G (x,t,z)$ and 
$\psi ([t,z]) = \widetilde \psi (t,z)$.
Then the perturbation triple~$(\psi,G,k)$ satisfies condition~(II)
of an equivariant extension. 
Moreover, since the perturbation pair~$(\chi,H)$ is of Moser type, 
the triple~$(\psi,G,k)$ also meets condition~(i) of admissibility.
It does not necessarily satisfy condition~(ii) of admissibility. 
However, we can remedy this by replacing $k$ by~$Ck$
for a large enough positive constant~$C$. 
This finishes the proof of the lemma. 
\proofend

\subsection{Proof of Theorem~A}

Assume that $\Sigma$ is displaceable in~$M$,
and choose a defining Hamiltonian~$F \colon M \to \RR$ for~$\Sigma$ meeting assumption~\eqref{e:ass}.
In view of the definition~\eqref{e:dirlim} of $\RFH^\TT (\Sigma, M)$, it suffices to show that
$\FH(\cA^{F,N;\TT}) =0$ for each~$N$.
So fix $N \in \NN$.

Choose $\chi \colon S^1 \to [0,\infty)$ with $\supp (\chi) \subset (0,\frac 12)$
and $\int_{S^1} \chi (t) \, dt =1$,
and choose a Hamiltonian function $H \colon M \times S^1 \to \RR$ with 
$H(\cdot,t) = 0$ for all $t \in [0, \frac 12]$ whose time 1-flow~$\varphi_H$ displaces~$\Sigma$.
By Lemma~\ref{l:extension}, the pair $(\chi,H)$ has an equivariant extension~$(\psi,G,k)$.
Let $[v,\eta, z]$ be a critical point of $\cA_{\psi,G,k}$.
Choose $z \in S^{2N+1}$ over $[z]$.
By Lemma~\ref{l:crit} and by~(II) of Definition~\ref{def.equiv.extension},
\begin{equation*} 
\left.
\begin{array}{rcl}
\dot v(t) &=& \eta \, \chi(t+t_z) \,X_F \bigl( v(t) \bigr) + X_{H(\cdot, t+t_z)} \bigl( v(t),t \bigr), \\ [0.4em]
0 &=& \int_{S^1} \chi(t+t_z)\, F \bigl(v(t) \bigr) \,dt .
\end{array}
\right\}
\end{equation*}
%
By Lemma~\ref{l:intersection}, $v(t_z)$ is a leafwise intersection point for $H(\cdot, t+t_z)$.
This is impossible because $\varphi_H$ displaces~$\Sigma$.
It follows that the functional $\cA_{\psi,G,k} = \cA^{F,N;\TT}_{\psi,G,k}$ has no critical points.
The Floer homology $\FH (\cA^{F,N;\TT}_{\psi,G,k})$ is defined along the lines of Section~\ref{ss:eqRFH},
see~Section~\ref{s:invariance}.
Since $\cA^{F,N;\TT}_{\psi,G,k}$ has no critical points, 
the Floer complex of $\cA^{F,N;\TT}_{\psi,G,k}$ is trivial, and hence $\FH (\cA^{F,N;\TT}_{\psi,G,k}) =0$.
Theorem~A thus follows from the invariance $\FH (\cA^{F,N;\TT}_{\psi,G,k}) \cong \FH(\cA^{F,N;\TT})$,
which is proven in the next section.

\section{Invariance}  \label{s:invariance}

The goal of this section is to prove

\begin{proposition} \label{p:invariance}
$\FH (\cA^{F,N;\TT}_{\psi,G,k}) \cong \FH(\cA^{F,N;\TT})$.
\end{proposition}

This isomorphism can be proven along the lines of the proof of Corollary~3.4 in~\cite{CieFra09}.
In this section we give a different proof.

We start with reviewing two continuation methods for showing invariance of a Floer-type homology.
For simplicity, we describe these methods in the setting of 
Morse homology and Morse--Bott homology on a non-compact manifold~$M$.
For $i =0,1$ let $f_i \colon M \to \RR$ be smooth Morse functions with compact
critical sets~$\Crit f_i$.

\m \ni
{\bf Method 1.}
Assume that there is a smooth family $\{ f_s \}_{0 \leqslant s \leqslant 1}$ of Morse functions $f_s \colon M \to \RR$ 
such that the critical sets $\Crit f_s$ are all isotopic.
More precisely, assume that there is a diffeomorphism
$$
\Psi \colon \Crit f_0 \times [0,1] \,\to\, \coprod_{0 \leqslant s \leqslant 1} \Crit f_s \times \{s\} ,
\qquad
\bigl( x, s \bigr) \mapsto \bigl( x_s, s \bigr) .
$$ 
For a Riemannian metric $j_s$ on~$M$ and for $x_s,y_s \in \Crit f_s$ denote by $\widehat \cM_{f_s,j_s}(x_s,y_s)$ the set of negative
gradient flow lines from~$x_s$ to~$y_s$, and by $\cM_{f_s,j_s}(x_s,y_s) := \widehat \cM_{f_s,j_s}(x_s,y_s)/\RR$ 
the space of unparametrized gradient flow lines.

\begin{figure}[h]
 \begin{center}
  \psfrag{0}{$0$}
  \psfrag{1}{$1$}
  \psfrag{s}{$s$}
  \psfrag{M}{$\cM_{f_s,j_s} \bigl( x_s, y_s \bigr)$}
  \leavevmode\epsfbox{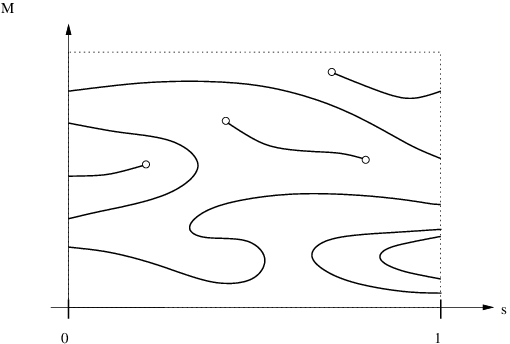}
 \end{center}
 \caption{The union of moduli spaces $\bigcup_{0 \leqslant s \leqslant 1} \cM_{f_s,j_s} \bigl( x_s, y_s \bigr) \times \{s\}$.}
 \label{fig.moduli}
\end{figure}

For $s=0,1$ choose a Riemannian metric~$j_s$ on~$M$ such that the pair $(f_s,j_s)$ is Morse--Smale.
Then one can define the Morse homology of~$f_s$ by counting elements of $\cM_{f_s,j_s}(x_s,y_s)$ 
for $x_s,y_s \in \Crit f_s$ with $\ind (x_s) = \ind (y_s)+1$, $s=0,1$.
For a generic smooth path of Riemannian metrics~$\{j_s\}$ from~$j_0$ to~$j_1$ 
and for $x,y \in \Crit f_0$ with $\ind (x) = \ind (y)+1$,
the union of moduli spaces
$$
\cM_{\{f,j\}} (\{x\},\{y\}) \,=\, \bigcup_{0 \leqslant s \leqslant 1} \cM_{f_s,j_s} \bigl( x_s, y_s \bigr) \times \{s\}
$$ 
is then a 1-dimensional smooth manifold with boundary that is ``transverse at 0 and 1'', 
i.e., for $s=0,1$ the points in $\cM_{f_s,j_s} \bigl( x_s, y_s \bigr) \times \{s\}$ belong to the boundary of 
$\cM_{\{f,j\}} (\{x\},\{y\})$, see Figure~\ref{fig.moduli}.
If one can show that the sets $\widehat \cM_{f_s,j_s} \bigl( x_s, y_s \bigr)$, $0 \leqslant s \leqslant 1$,
are uniformly bounded,  
it follows that the Morse homologies of~$f_0$ and of~$f_1$ are isomorphic.

%
%

Indeed, $\cM_{\{f,j\}} (\{x\},\{y\})$ is the union $\cM_1 \coprod \cM_2$ of two types of components:
The components of~$\cM_1$ are compact intervals with boundary over~$0$ and~$1$,
and the components of~$\cM_2$ are half-open intervals (with boundary over $0$ or~$1$) 
or open intervals.
If~$\cM_2$ is empty, then the coefficients $\nu (x_i,y_i^k) = \# \cM_{f_i,j_i} (x_i,y_i^k) \mod 2$
in the boundary operator 
$$
\partial_i \2 x_i \,=\, \sum_k \nu (x_i,y_i^k) \,y_i^k
$$
are the same for $i=0,1$.  
The components of~$\cM_2$ may change the coefficients $\nu (x_i,y_i^k)$,
but they do not alter the Morse homology.
Indeed, the contribution of the components of~$\cM_2$  to the boundary operator can be computed explicitely, 
and from this one can write down an explicit chain homotopy equivalence between the Morse chain complexes of
$(f_0,j_0)$ and $(f_1,j_1)$, see \cite[Lemmata 3.5 and~3.6]{Flo88:Lag}.
We illustrate this by an example:

Suppose $\Crit f_s$ has three critical points, $a_s,b_s$ of index~$1$ and $c_s$ of index~$0$.
Suppose that at $s=0$ there is exactly one gradient flow line~$\gamma_0$, from $b_0$ to~$c_0$.
Then the Morse homology is generated by~$a_0$: 
$$
\MH (f_0,j_0) \,=\, \MH_1(f_0,j_0) = \ZZ_2 \langle a_0 \rangle .
$$
Assume now that at some time $s^* \in (0,1)$ a gradient flow line $\gamma_{ab}$ from $a_{s^*}$ to~$b_{s^*}$ appears.
This flow line is not generic, and immediately disappears.
The flow line $\gamma_{ab}$ affects the two families of moduli spaces 
$\cM_{f_s,j_s}(b_s,c_s)$ and $\cM_{f_s,j_s}(a_s,c_s)$ as follows:
The moduli spaces $\cM_{f_s,j_s}(b_s,c_s)$ are not affected:
Before time~$s^*$ this space contains exactly one gradient flow line~$\gamma_s$, which persists beyond time~$s^*$.

\begin{figure}[h]
 \begin{center}
  \psfrag{a0}{$a_0$}
  \psfrag{b0}{$b_0$}
  \psfrag{c0}{$c_0$}
  \psfrag{a1}{$a_1$}
  \psfrag{b1}{$b_1$}
  \psfrag{c1}{$c_1$}
  \psfrag{as}{$a_{s^*}$}
  \psfrag{bs}{$b_{s^*}$}
  \psfrag{cs}{$c_{s^*}$}
  \psfrag{g0}{$\gamma_0$}
  \psfrag{g1}{$\gamma_1$}
  \psfrag{gs}{$\gamma_{s^*}$}
  \psfrag{ga}{$\gamma_{ab}$}
  \leavevmode\epsfbox{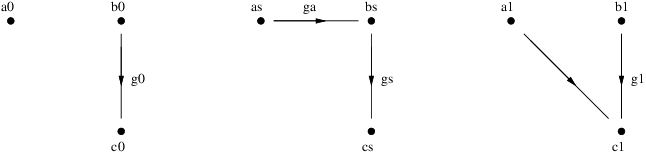}
 \end{center}
 \caption{The gradient flow lines at $s=0$, $s=s^*$, $s=1$.}
 \label{fig.exlines}
\end{figure}
%
%

The moduli spaces $\cM_{f_s,j_s}(a_s,c_s)$ were empty for $s < s^*$.
At time $s^*$ there is a broken gradient flow line from $a_{s^*}$ to~$c_{s^*}$,
namely $\gamma_{ab}$ followed by~the gradient flow line~$\gamma_{s^*}$ from $b_{s^*}$ to $c_{s^*}$.
These two flow lines can be glued together to a unique gradient flow line from $a_s$ to~$c_s$.
Hence $\nu (a_s, c_s)$ changes at~$s^*$ from~$0$ to~$1$.
For $s>s^*$ we now have one gradient flow line from $a_s$ to~$c_s$ and one from~$b_s$ to~$c_s$.
But this change does not affect the Morse homology: $c_s$ is still in the image of the boundary operator~$\pp_s$,
and while now neither $a_s$ nor~$b_s$ are in the kernel, $a_s-b_s$ is in the kernel of~$\pp_s$.
Hence we still have 
$$
\MH (f_1,j_1) \,=\, \MH_1(f_1,j_1) = \ZZ_2 \langle a_1-b_1 \rangle .
$$

\begin{figure}[h]
 \begin{center}
  \psfrag{Mbc}{$\cM_{f_s,j_s}(b_s,c_s)$}
  \psfrag{Mac}{$\cM_{f_s,j_s}(a_s,c_s)$}
  \psfrag{0}{$0$}
  \psfrag{1}{$1$}
  \psfrag{s}{$s^*$}
  \psfrag{gs}{$\gamma_{s^*}$}
  \leavevmode\epsfbox{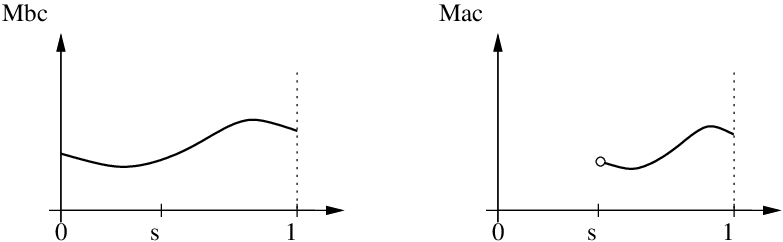}
 \end{center}
 \label{fig.exmod}
\end{figure}
%
%

A bifurcation as above, that creates a component in~$\cM_2$, 
is called a slide bifurcation, or a handle slide, 
since such a bifurcation acts on the corresponding handle decomposition of~$M$ 
by sliding one handle over another.
The other type of bifurcation that appears in a generic isotopy between Morse functions
are birth bifurcations and death bifurcations, namely the birth of two critical points
or the cancellation of two critical points.
Such bifurcations do not arise in the situation at hand.

\m 
Below we shall apply this method in a Morse--Bott set-up:
Assume there is a smooth family $\{f_s\}_{0 \leqslant s \leqslant 1}$ of Morse--Bott functions~$f_s \colon M \to \RR$ 
with compact critical sets~$\Crit f_s$
and a diffeomorphism
$$
\Psi \colon \Crit f_0 \times [0,1] \,\to\, \coprod_{0 \leqslant s \leqslant 1} \Crit f_s \times \{s\} ,
\qquad
\bigl( x, s \bigr) \mapsto \bigl( x_s, s \bigr) .
$$ 
Choose a Morse function $h_0$ on~$\Crit f_0$. Then the functions
$$
h_s(x_s) \,:=\, h_0(x)
$$
are Morse functions on~$\Crit f_s$, and the sets $\Crit h_s$ are isotopic.
For a Riemannian metric~$g_s$ on~$\Crit f_s$,  
for a Riemannian metric $j_s$ on~$M$ and for $x_s,y_s \in \Crit h_s$ denote by $\widehat \cM_{f_s,j_s,h_s,g_s}(x_s,y_s)$ 
the set of negative gradient flow lines with cascades from~$x_s$ to~$y_s$, and by 
$\cM_{f_s,j_s,h_s,g_s}(x_s,y_s)$ the space of unparametrized gradient flow lines with cascades.

For $s=0,1$ choose a Riemannian metric~$g_s$ on~$\Crit f_s$
such that the pair $(h_s,g_s)$ is Morse--Smale.
For generic Riemannian metrics~$j_s$ on~$M$ one can then define the Morse--Bott homology of the quadruples~$(f_s,j_s,h_s,g_s)$, 
$s=0,1$, by counting elements of the 0-dimensional components $\cM_{f_s,j_s,h_s,g_s}(x_s,y_s)$, 
see~\cite[Appendix~A]{Fra04}.
For a generic smooth path of Riemannian metrics~$\{g_s\}$ on~$\Crit f_s$ and 
for a generic smooth path of Riemannian metrics~$\{j_s\}$ on~$M$
from~$(j_0,g_0)$ to~$(j_1,g_1)$,
we have that 
for each pair $x,y \in \Crit h_0$ for which $\cM_{f_0,j_0,h_0,g_0}(x,y)$ is $0$-dimensional,
the union of moduli spaces
$$
\cM_{\{f,j,h,g\}} (\{x\},\{y\}) \,=\, \left\{ (u,s) \mid u \in \cM_{f_s,j_s,h_s,g_s} \bigl( x_s, y_s \bigr) ; \; 0 \leqslant s \leqslant 1 \right\} ,
$$ 
is a 1-dimensional smooth manifold with boundary that is ``transverse at 0 and 1''.
If one can show that the sets $\widehat \cM_{f_s,j_s,g_s,h_s} \bigl( x_s, y_s \bigr)$, $0 \leqslant s \leqslant 1$,
are uniformly bounded,  
it follows that the Morse homologies of~$f_0$ and of~$f_1$ are isomorphic.

\m \ni
{\bf Method 2 (Floer continuation).}
Choose a smooth monotone function $\beta \colon \RR \to [0,1]$ with $\beta (s) = 0$ for $s \leqslant 0$
and $\beta (s) = 1$ for $s \geqslant 1$.
For $s \in \RR$ define the function 
$$
f_s \,=\, (1-\beta(s)) f_0 + \beta(s) f_1 .
$$
For $x \in \Crit f_0$ and $y \in \Crit f_1$ and 
for a smooth family of Riemannian metrics $\{g_s\}$ with $g_s = g_0$ for $s \leqslant 0$ and $g_s = g_1$ for $s \geqslant 1$
consider the gradient equation with asymptotic boundary conditions
\begin{equation} \label{e:Morse}
\left\{ 
\begin{array}{l}
\dot u (s) = -\nabla_{g_s} f_s (u(s)), \quad s \in \RR;  \\ [0.4em]
\displaystyle\lim_{s \to -\infty} u(s) = x, \quad \lim_{s \to \infty} u(s) = y.
\end{array}
\right.
\end{equation}
For a generic choice of the path $\{g_s\}$ and for $x \in \Crit f_0$ and $y \in \Crit f_1$ 
with $\ind (x) = \ind (y)$,
the space of solutions to~\eqref{e:Morse} is a smooth 0-dimensional manifold.
If one can show that this space is bounded, then it is finite.
Counting these solutions then defines a chain homomorphism between the Morse chain complexes of~$f_0$ and~$f_1$, 
that induces an isomorphism between the Morse homologies of~$f_0$ and~$f_1$.

\smallskip
Similarly, given triples $(j_s, h_s, g_s)$ for $s=0,1$ with $(h_s,g_s)$ Morse--Smale pairs and $j_s$ generic, 
Floer continuation can be used to show that the Morse homologies of~$(f_0,j_0, h_0, g_0)$ and~$(f_1,j_1, h_1, g_1)$
are isomorphic, see~\cite[Theorem~A.17]{Fra04}.

\m \ni 
{\bf Historical Remark.}
{\rm
Floer used Method~1 in~\cite{Flo88:Lag} to prove invariance of his homology for Lagrangian intersections.
(He also dealt with isolated bifurcations of the critical sets, namely birth and death bifurcations, 
by first putting them into normal form and then constructing a chain map between the complex before and after
the bifurcation that induces an isomorphism in homology.) 
Such a bifurcation analysis was later also used in~\cite{EES05, Lee05, Sul02}.)
The powerful and flexible Method~2 was invented by Floer only later in~\cite{Flo89:hol}.
\diam
}

\m
Proposition~\ref{p:invariance} can be proven by Method 2, by adapting the proof of Corollary~3.4 in~\cite{CieFra09}.
We leave the minor modifications to the interested reader. 
Here we give a different argument, that takes into account the structure of the functional~$\cA^{F,N;\TT}_{\psi,G,k}$,
and uses Method~1 once and Method~2 twice. 

\m
Consider the four functionals on~$\cL \times \RR \times_{\TT} S^{2N+1}$,
\begin{eqnarray*}
\cA_0 ([v,\eta,z]) &=& -\int_{S^1} v^* \lambda 
-\eta \int_{S^1}  F \bigl( v(t) \bigr) \,dt , \phantom{\widetilde \psi \left(t,z \right) \,
-\int_{S^1} \widetilde G \bigl( v(t),t,z \bigr) \, dt - \widetilde k(z) } \\ [0.2em]
\cA_1 ([v,\eta,z]) &=& -\int_{S^1} v^* \lambda 
-\eta \int_{S^1}  F \bigl( v(t) \bigr) \,dt \phantom{\widetilde \psi \left(t,z \right)\,}
\phantom{-\int_{S^1} \widetilde G \bigl( v(t),t,z \bigr) \, dt} \, - \widetilde k(z) , \\ [0.2em]
\cA_2 ([v,\eta,z]) &=& -\int_{S^1} v^*\lambda 
-\eta \int_{S^1} \widetilde \psi \left(t,z \right) F \bigl( v(t) \bigr) \,dt 
\phantom{-\int_{S^1} \widetilde G \bigl( v(t),t,z \bigr) \, dt} \, - \widetilde k(z) , \\ [0.2em]
\cA_3 ([v,\eta,z]) &=&  
-\int_{S^1} v^*\lambda 
-\eta \int_{S^1} \widetilde \psi \left(t,z \right) F \bigl( v(t) \bigr) \,dt 
-\int_{S^1} \widetilde G \bigl( v(t),t,z \bigr) \, dt - \widetilde k(z) .
\end{eqnarray*}
The functionals $\cA_0$ and $\cA_1$ are Morse--Bott by our assumption~\eqref{e:ass} on~$\Sigma$
and since $k$ is Morse, 
while $\cA_2$ is Morse--Bott by Lemma~\ref{l:MB} below.
The functional~$\cA_3$ is Morse--Bott because it has no critical points.
Hence the four lifted functionals $\widetilde \cA_i \colon \cL \times \RR \times S^{2N+1} \to \RR$ 
are also Morse--Bott.

The Floer homology $\FH(\cA_0) = \FH(\cA_0,h_0,J_0)$ was defined in Section~\ref{ss:eqRFH},
and the Floer homology~$\FH(\cA_i)$ for $i =1,2,3$ is defined in the same way:
One chooses a Morse function~$h_i$ and a Riemannian metric~$g_i$ on $\Crit \cA_i$
such that $(h_i,g_i)$ is a Morse--Smale pair, 
lifts them to the Morse--Bott function~$\widetilde h_i$ and the $\TT$-invariant metric~$\widetilde g_i$
on~$\Crit \widetilde \cA_i$, and defines the boundary of a critical point~$c^+$ of~$h_i$
by counting rigid $\TT$-families of unparametrized negative gradient flow lines with cascades in~$\cM (c^+,c^-)$ 
between critical $\TT$-orbits~$C^+$ and~$C^-$ of~$\widetilde h_i$,
with respect to the $\TT$-invariant Riemannian metric~$\widetilde g_i$ on~$\Crit \widetilde \cA_i$
and a generic family~$J_{t,z}(\cdot, \eta)$ in~$\mathcal{J}^{S^1}$.

It follows from Method~2 that $\FH(\cA_0) \cong \FH(\cA_1)$ and that $\FH(\cA_2) \cong \FH(\cA_3)$.
This is easy for the passage $\cA_0 \leadsto \cA_1$: 
The summand $\widetilde k(z)$ is bounded with all its derivatives.
The claim thus follows from the $L^\infty$-bound on each space $\widehat \cM (c^+,c^-)$ 
of gradient flow lines with cascades between a pair of critical circles of~$\widetilde h_0$ 
given in the proof of Corollary~3.3 in~\cite{CieFra09}.
For the passage $\cA_2 \leadsto \cA_3$, invariance follows as in \cite[Section 2]{AlbFra10}, 
by either choosing~$\widetilde G$ sufficiently small in~$L^\infty$ (which we are free to do)
or by decomposing the isotopy $\cA_2 \leadsto \cA_3$ into many small isotopies.

The isomorphism $\FH(\cA_1) \cong \FH(\cA_2)$ can also be shown by applying
Method~2 to the parts of a sufficiently fine decomposition of the isotopy $\cA_1 \leadsto \cA_2$ 
(see the proof of Corollary~3.4 in~\cite{CieFra09}).
This argument is somewhat harder, since $\eta$ appears in front of the summand that is altered. 
We circumvent this difficulty by applying Method~1.
Choose a smooth monotone function $\beta \colon [1,2] \to [0,1]$ 
with $\beta (s)=0$ for $s$ near~$1$ and $\beta (s)=1$ for $s$ near~$2$.
For $s \in [1,2]$ set 
$$
\widetilde \psi_s (t,z) \,=\, \bigl( 1-\beta (s) \bigr) \cdot 1 + \beta(s) \cdot \widetilde \psi (t,z) \,=\, 
1+\beta(s) \bigr( \widetilde \psi (t,z)-1 \bigl) .
$$
Then $\int_{S^1} \widetilde \psi_s(t,z)\, dt =1$ for all~$s$.
Consider the family of functionals
$$
\widetilde \cA_s (v,\eta,z) \,:=\, -\int_{S^1} v^*\lambda - \eta \int_{S^1} \widetilde \psi_s (t,z)\, F(v(t))\, dt + \widetilde k(z),
\qquad 1 \leqslant s \leqslant 2 .
$$
Then $\widetilde\cA_s = \widetilde\cA_1$ for $s$ near~$1$ and $\widetilde\cA_s = \widetilde\cA_2$ for $s$ near~$2$.

The critical manifolds $\Crit \widetilde\cA_s$ are in canonical bijection with $\Crit \widetilde\cA_1$.
Indeed, looking at~\eqref{crit.tilde} with $\widetilde G = 0$ and $\widetilde \psi$ replaced by~$\widetilde \psi_s$, 
we see that they all contain $\Sigma \times \{0\} \times \Crit \widetilde k$.
Moreover, given $z \in \Crit \widetilde k$, and with 
$$
s_z(t) \,:=\, \int_0^t \widetilde \psi_s(\tau, z)\, d \tau ,
$$
the periodic orbit $(v(t), \eta, z)$ of~$X_F$ with period~$|\eta|$ corresponds to the reparametrized 
orbit $(v(s_z(t)), \eta, z)$ of~$\widetilde \psi_s(t, z)\,X_F$ with period~$|\eta|$. 
(The orbit $v(s_z(t))$ also has period~$|\eta|$ because $s_z(1)=1$.)
More formally, the reparametrization map 
$$
\widetilde \Psi \colon
\Crit \widetilde \cA_1 \times [1,2] \,\to\, \coprod_{1 \leqslant s \leqslant 2} \Crit \widetilde \cA_s \times \{s\}, 
\qquad
\bigl( (v(\cdot), \eta, z), s \bigr) \mapsto \bigl( (v(s_z(\cdot)), \eta, z), s \bigr)
$$ 
is a diffeomorphism.

\begin{lemma} \label{l:MB}
For each $s \in [1,2]$ the critical set $\Crit \widetilde \cA_s$ is a Morse--Bott submanifold of~$\widetilde \cA_s$.
\end{lemma}

Before giving the proof, we use the lemma to prove Proposition~\ref{p:invariance}.
All the functionals~$\widetilde \cA_s$ and all the sets $\Crit \widetilde \cA_s$ are $\TT$-invariant.
Choose a Morse function $h_1$ on~$\Crit \cA_1$. Then the functions
$$
h_s \left( [v(s_z(\cdot)),\eta,z] \right) \,:=\, h_1([v(\cdot), \eta,z])
$$
are Morse functions on~$\Crit \cA_s$, and the sets $\Crit h_s$ are isotopic.

For a Riemannian metric~$g_s$ on~$\Crit \cA_s$,  
for a family $\bJ_s := (J_{t,z}(\cdot, \eta))_s$ in~$\mathcal{J}^{S^1}$ and for $c_s^+,c_s^- \in \Crit h_s$,
denote by $\widehat \cM_{\cA_s,\bJ_s,h_s,g_s}(c_s^+,c_s^-)$ 
the set of negative gradient flow lines with cascades from~$c_s^+$ to~$c_s^-$, and by 
$\cM_{\cA_s,\bJ_s,h_s,g_s}(c_s^+,c_s^-)$ the space of unparametrized $\TT$-families of gradient flow lines with cascades,
as constructed in~Section~\ref{ss:eqRFH}.3.

For $s=1,2$ choose $g_s$ and~$\bJ_s$ as in the definition of the Floer homologies~$\FH(\cA_s)$:
$g_s$ is a Riemannian metric on~$\Crit \cA_s$ such that $(h_s, g_s)$ is a Morse--Smale pair,
and $\bJ_s$ is a generic family in~$\mathcal{J}^{S^1}$.
Then for a generic smooth path of Riemannian metrics $\{g_s\}$ on~$\Crit \cA_s$
and for a generic smooth path of families $\{\bJ_s\}$ in~$\mathcal{J}^{S^1}$ 
from $(g_1, \bJ_1)$ to~$(g_2,\bJ_2)$,
we have that for each pair $c^+,c^- \in \Crit h_1$ for which $\cM_{\cA_1,\bJ_1,h_1,g_1}(c^+,c^-)$ is $0$-dimensional,
the union of moduli spaces
$$
\cM_{\cA_,\bJ,h,g} (\{c^+\},\{c^-\}) \,=\, 
\bigcup_{1 \leqslant s \leqslant 2} \cM_{\cA_s,\bJ_s,h_s,g_s} \bigl( c_s^+,c_s^- \bigr)  \times \{s\}
$$ 
is a 1-dimensional smooth manifold that is ``transverse at 0 and 1''.
Notice that the map $\widetilde \Psi$ is action-preserving: 
$\widetilde \cA_s (x_s) = \widetilde \cA_1 (x)$.
The space $\widehat \cM_{\cA_1,\bJ_1,h_1,g_1} (c^+,c^-)$ is $L^\infty$-bounded,
and in fact there is a uniform $L^\infty$-bound on the spaces $\widehat \cM_{\cA_s,\bJ_s,h_s,g_s} \bigl( c^+_s, c^-_s \bigr)$,
$1 \leqslant s \leqslant 2$, see the proof of Corollary~3.3 in~\cite{CieFra09}.
It follows that $\FH(\cA_1) \cong \FH(\cA_2)$.

\m
\ni
{\it Proof of Lemma~\ref{l:MB}.}
We use the method in Appendix~A.1 of~\cite{AlbFra10}.
Fix $s$, and fix a critical point $(v_0, \eta_0,z_0) \in \cL \times \RR \times S^{2N+1}$.
We decompose $\widetilde \cA_s$ as
$$
\widetilde \cA_s (v,\eta,z) \,=\, \cA_0 (v)  + \eta_0 \, \cF_{\Delta \psi} (v,z) + (\eta_0-\eta) \2 \cF_\psi (v,z) + \widetilde k (z)
$$
where 
\begin{eqnarray*}
\cA_0 (v)         &:=&  -\int_{S^1} v^* \lambda - \eta_0 \int_{S^1} \widetilde \psi_s (t,z_0) F(v(t)) \, dt , \\ [0.2em]
\cF_{\Delta \psi} (v,z)  &:=& \int_{S^1} \left( \widetilde \psi_s (t,z_0)-\widetilde \psi_s (t,z)\right) F(v(t)) \, dt , \\ [0.2em]
\cF_{\psi} (v,z)         &:=& \int_{S^1} \widetilde \psi_s (t,z) F(v(t)) \, dt .
\end{eqnarray*}
In order to compute the Hessian of~$\widetilde \cA_s$ at~$(v_0, \eta_0,z_0)$, we apply ``a change of coordinates'':
Consider the twisted loop space 
$$
\cL_{\eta_0F} \,:=\, \left\{ w \in C^\infty([0,1], M) \mid w(0) = \phi^1_{\eta_0F}(w(1)) \right\}
$$
and the diffeomorphism $\Phi_{\eta_0F} \colon \cL_{\eta_0F} \to \cL = C^\infty(S^1,M)$ given by
$$
\Phi_{\eta_0F}(w)(t) \,=\, \phi_{\eta_0 F_s}^t (w(t))
$$
where we abbreviated $F_s (\cdot) := \widetilde \psi_s (t,z_0) F(\cdot)$.
Then the path $w_0 = \Phi^{-1}_{\eta_0F} \circ v_0 = v_0(0) \in \Sigma$ is constant.
Hence tangent vectors $\hat w(t)$ at $w_0$ are curves in the linear space $T_{w_0}M$ with
\begin{equation} \label{e:twist.linear}
\hat w (1) \,=\, d \phi_{\eta_0F_s}^{-1}(w_0) \,\hat w(0) .
\end{equation}

We are going to compute the kernel of the Hessian of the pulled-back functional 
$$
\cA_s^\Phi := (\Phi_{\eta_0F} \times \id_\RR \times \id_{S^{2N+1}})^* \widetilde \cA_s 
\colon \, \cL_{\eta_0F} \times \RR \times S^{2N+1} \to \RR
$$
at the critical point $(w_0, \eta_0,z_0)$.
First notice that
$\Phi_{\eta_0F}^* d\cA_0 (w)[\hat w] = \int_0^1 \omega (\frac{d}{dt}w, \hat w) dt$
for any $w \in \cL_{\eta_0F}$ and $\hat w \in T_w \cL_{\eta_0F}$, 
and that
\begin{eqnarray*}
(\Phi_{\eta_0F} \times \id_{S^{2N+1}})^* \cF_{\Delta \psi} (w,z) &=& 
\int_0^1 \left( \widetilde \psi_s (t,z_0) - \widetilde \psi_s(t,z) \right) F(w(t)) \, dt , \\
(\Phi_{\eta_0F} \times \id_{S^{2N+1}})^* \cF_\psi (w,z) &=& 
\int_0^1 \widetilde \psi_s(t,z) \, F(w(t)) \, dt 
\end{eqnarray*}
(since $F$ is preserved under $\phi_{\eta_0F_s}^t$).
The differential of $\cA_s^\Phi$ therefore is
\begin{eqnarray*}
d \cA_s^\Phi (w,\eta,z)[\hat w, \hat \eta,  \hat z] &=&
\int_0^1 \omega (\tfrac{d}{dt} w, \hat w)\, dt \\
&& + \eta_0 \int_0^1 \left\{ \left( \widetilde \psi_s (t,z_0) - \widetilde \psi_s (t,z) \right) dF(w(t)) \, \hat w(t) 
                                   -\pp_z \widetilde \psi_s (t,z) \,\hat z \, F(w(t)) \right\} dt \\
&& - \hat \eta \int_0^1 \widetilde \psi_s (t,z) \,F(w(t)) \,dt \\
&& + (\eta_0-\eta) \int_0^1 	\widetilde \psi_s (t,z) \,dF(w(t)) \, \hat w(t) + 
                                               \pp_z \widetilde \psi_s (t,z) \,\hat z \, F(w(t)) \, dt \\			
&& + d \widetilde k (z) \hat z .
\end{eqnarray*}
At the critical point $x_0 := (w_0,\eta_0,z_0)$ 
the Hessian of $\cA_s^\Phi$ applied to $\xi_i := (\hat w_i, \hat \eta_i, \hat z_i)$ 
therefore is
\begin{eqnarray*}
\Hess \, \cA_s^\Phi (x_0)[\xi_1, \xi_2] &=&
\int_0^1 \omega (\tfrac{d}{dt} \hat w_1, \hat w_2)\, dt \\
&& - \eta_0 \int_0^1 \left\{ \pp_z \widetilde \psi_s (t,z_0) \, \hat z_1 \, dF(w_0) \, \hat w_2(t) 
                           + \pp_z \widetilde \psi_s (t,z_0) \, \hat z_2 \, dF(w_0) \, \hat w_1 (t) \right\} dt \\
&& - \hat \eta_1 \int_0^1 \widetilde \psi_s (t,z_0) \,dF(w_0) \, \hat w_2(t) \,dt \\
&& - \hat \eta_2 \int_0^1 \widetilde \psi_s (t,z_0) \,dF(w_0) \, \hat w_1(t) \, dt  \\			
&& + \Hess \, \widetilde k (z) (\hat z_1, \hat z_2) 
\end{eqnarray*}
where we have used that $F(w_0)=0$.
A tangent vector $(\hat w, \hat \eta, \hat z)$ therefore belongs to the kernel of $\Hess \, \cA_s^\Phi (w_0,\eta_0,z_0)$
if and only if
\begin{eqnarray} 
0 &=& 
\tfrac{d}{dt} \hat w(t) - \hat \eta \, \widetilde\psi_s (t,z_0) \, X_F (w_0) 
                       - \eta_0 \, \partial_z \widetilde \psi_s(t,z_0) \, \hat z \, X_F(w_0) ,
\label{Hess1} \\  [.4em]
0 &=& \int_0^1 \widetilde \psi_s(t,z_0)\, dF(w_0) \2 \hat w(t) \,dt , 
\label{Hess2} \\ 
0 &=& -\eta_0 \int_0^1  dF(w_0) \2 \hat w(t) \, \pp_z \widetilde \psi_s (t,z_0) (\cdot) \,dt 
                          +  \Hess\2 \widetilde k(z_0) (\hat z, \cdot) .
\label{Hess3}
\end{eqnarray}
Denote by $H_{z_0} = \left\{ \tau z_0 \mid \tau \in S^1 \right\}$ the Hopf circle in~$S^{2N+1}$ through~$z_0$. 

Assume first that $\eta_0=0$. Then $(v_0,\eta_0,z_0)$ belongs to the critical component $\Sigma \times \{0\} \times H_{z_0}$ of 
``constant in~$\Sigma$ loops''.
Since $\eta_0 =0$, \eqref{Hess3} yields $\hat z \in T_{z_0}H_{z_0}$,
and integrating~\eqref{Hess1} yields
$$
\hat w (1) \,=\, \hat w(0) + \hat \eta \, X_F(w_0)
$$
(since $s_{z_0}(1)=1$).
Since in this case $\Phi_{\eta_0F} \colon \cL \to \cL$ is the identity mapping, $\hat w(1)=\hat w(0)$, and so $\hat \eta =0$.
By now, \eqref{Hess1} reads $\tfrac{d}{dt} \hat w(t) =0$, that is, $\hat w(t) \equiv \hat w(0) \in T_{w_0}M$ is constant.
Finally, \eqref{Hess2} shows that $\hat w(0) \in T_{w_0} \Sigma$.
The kernel of the Hessian of $\widetilde \cA_s$ at $(v_0,\eta_0,z_0) = (w_0,0,z_0)$ is thus identified with 
$T_{w_0}\Sigma \times T_{z_0}H_{z_0}$.

Assume now that $\eta_0 \neq 0$.
Then $S_{v_0} := \left\{ v_0(\cdot -\tau) \mid \tau \in S^1 \right\}$ is an embedded circle in~$\cL$.
Hence the critical component of~$(v_0, \eta_0, z_0)$ is the torus $S_{v_0} \times \{\eta_0\} \times H_{z_0}$.
It is clear that the kernel of the Hessian of $\widetilde \cA_s$ at $(v_0,\eta_0,z_0)$ has dimension at least two, 
and we must show that the dimension is two.
By assumption~\eqref{e:ass}, $1$ has multiplicity~$2$ in the spectrum of $d \phi_{\eta_0 F}^{-1}(w_0)$.
Since $\phi_{\eta_0 F_s} = \phi_{\eta_0 F}$, the same holds true for 
$L_s := d \phi_{\eta_0 F_s}^{-1}(w_0)$.
Recall that $s_z(t) = \int_0^t \widetilde \psi_s(\tau, z)\, d\tau$.
Integrating \eqref{Hess1} we get
\begin{equation} \label{e:wt}
\hat w (t) \,=\, \hat w(0) + \hat \eta \, s_{z_0}(t) \2 X_F(w_0) + \eta_0  \, \pp_z \big|_{z_0} s_z(t) \, \hat z \, X_F(w_0) .
\end{equation}
In particular (since $s_z(1) = 1$ for all~$z$), and by~\eqref{e:twist.linear},
\begin{equation} \label{e:w1}
\hat w(1) \,=\, \hat w(0) + \hat \eta \, X_F (w_0) \,=\, L_s\, \hat w(0) .
\end{equation} 
Consider the sub-vector space $V$ of~$T_{w_0}M$ spanned by $\hat w(0)$ and $X_F(w_0)$.
Assume that $V$ is 2-dimensional. Then \eqref{e:w1} and the fact that $1$ has multiplicity~$2$ 
in the spectrum of~$L_s$ show that $V$ is the whole 1-eigenspace of~$L_s$.
In particular, $V$ is symplectic.
On the other hand, since $dF (X_F)=0$, equations~\eqref{e:wt} and~\eqref{Hess2} show that
$$
dF(w_0) \, \hat w(0) \,=\, \int_0^1 \widetilde \psi_s (t,z_0)\, dF(w_0) \, \hat w(0)\, dt \,=\, 
                           \int_0^1 \widetilde \psi_s (t,z_0)\, dF(w_0) \, \hat w(t)\, dt \,=\, 0
$$
and hence $\hat w(0) \in T_{w_0} \Sigma$.
Since $X_F(w_0)$ generates the kernel of $\omega |_{T_{w_0} \Sigma}$, this contradicts $V$ being symplectic.

It follows that $\hat w(0) = r \1X_F(w_0)$ for some $r \in \RR$.
In particular, $L_s \,\hat w(0) = \hat w(0)$.
The second equation in~\eqref{e:w1} thus shows that $\hat \eta =0$.
Since $\hat w (0) \in V = \spa \1 (X_F(w_0))$,
equation~\eqref{e:wt} shows that $\hat w(t) \in V$ for all~$t$.
Therefore \eqref{Hess3} gives $\hat z \in \ker \Hess \, \widetilde k(z_0) = T_{z_0}H_{z_0}$.
We conclude with~\eqref{e:wt} that the kernel of~$\Hess\, \cA_s^\Phi (w_0,\eta_0,z_0)$ is
$$
\left\{ \left( \hat w(t), 0, \hat z \right) \mid \hat z \in T_{z_0}H_{z_0} \right\} =
\left\{ \bigl( r + \eta_0  \, \pp_z \big|_{z_0} s_z(t) \, \hat z \bigr) \, X_F(w_0), 0, \hat z \bigr) 
\mid r \in \RR,\, \hat z \in T_{z_0}H_{z_0} \right\} .
$$
Hence $\dim \ker \Hess\, \cA_s^\Phi (w_0,\eta_0,z_0) = \dim \ker \Hess\, \widetilde \cA_s (v_0, \eta_0,z_0) =2$.
\proofend


\section{Other approaches} \label{s:other}

In this note we have defined $\TT$-equivariant Rabinowitz--Floer homology $\RFH^\TT(\Sigma,M)$
via the Borel construction and Floer homology with cascades,
and we have proven the vanishing of $\RFH^\TT(\Sigma,M)$ for displaceable~$\Sigma$ by a leave-wise intersection argument.
There are several other approaches to construct a $\TT$-equivariant Rabinowitz--Floer homology
(two are mentioned in Remark~\ref{rem:grading}, and one more is outlined in 3.\ below),
all of which are expected to give the same result.
And there are different ways to prove the vanishing of $\RFH^\TT(\Sigma,M)$ or of these other versions 
for displaceable~$\Sigma$.
In particular, the arguments in 1.\ and~2.\ below imply the vanishing of the version 
defined in~\cite{CieOan15}, see~items (4) and~(3) on page~70 of~\cite{CieOan15}.
Let $V$ be the bounded component of $M \setminus \Sigma$,
and denote by~$\SH_*(V)$ its symplectic homology and  
by~$\SH_*^\TT(V)$ its equivariant symplectic homology.

\m
\ni
{\bf 1. Vanishing of $\RFH^\TT(\Sigma,M)$ via vanishing of~$\SH^\TT(V)$.}
There should be a~$\TT$-equivariant version of the long exact sequence
$$
\cdots \,\longrightarrow\, \SH^{-*}(V) \,\longrightarrow\, \SH_*(V) \,\longrightarrow\, \RFH_*(\Sigma,M) 
\,\longrightarrow\, \SH^{-*+1}(V) \,\rightarrow\, \cdots
$$
from~\cite{CieFraOan10}.
The vanishing of $\RFH^\TT(\Sigma,M)$ for displaceable~$\Sigma$ 
would then follow from the vanishing of~$\SH^\TT(V)$ proven in~\cite{BouOan12}.

\m
\ni
{\bf 2. Vanishing of $\RFH^\TT(\Sigma,M)$ via vanishing of $\RFH(\Sigma,M)$.}
It is shown in~\cite[Theorem~1.2]{BouOan12} that 
$$
\SH(V) =0 \,\Longleftrightarrow\, \SH^\TT(V)=0.
$$
While the implication $\Longleftarrow$ follows from the Gysin exact sequence 
in~\cite{BouOan09:Gysin},
the implication $\Longrightarrow$ follows from the fact that $\SH^\TT(V)$ is the limit 
of a spectral sequence whose second page is the tensor product of the homology 
of the classifying space~$BS^1$ and of~$\SH(V)$, \cite[\S 2.2]{BouOan12}.
It is expected that these two algebraic constructions can be adapted 
to Rabinowitz--Floer homology (cf.\ \cite[p.\ 6]{BouOan12}). 
Then 
$$
\RFH(\Sigma,M) =0 \,\Longleftrightarrow\, \RFH^\TT(\Sigma,M) =0.
$$
In particular, the vanishing of $\RFH^\TT(\Sigma,M)$ for displaceable~$\Sigma$ would then follow 
from the vanishing of~$\RFH(\Sigma,M)$ proven in~\cite{CieFra09}.
Together with the equivalence from~\cite[Theorem~13.3]{Rit13} we could then conclude
the equivalences
$$
\RFH^\TT(\Sigma,M) =0 \,\Longleftrightarrow\, \RFH(\Sigma,M) =0 
                      \,\Longleftrightarrow\, \SH(V) =0 
                      \,\Longleftrightarrow\, \SH^\TT(V) =0.
$$

\m
\ni
{\bf 3. Chekanov's construction of $S^1$-equivariant Floer homology.}
In the Borel-construction, approximations $S^{2N+1}$ of the classifying space $S^\infty = ES^1$
are somewhat clumsily added to the loop space as direct summands.
In Chekanov's version of $S^1$-equivariant Floer homology,
$S^{2N+1}$ does not appear as a space, but is incorporated into the boundary operator:
In the setting of Morse theory for a function $f \colon M \to \RR$ on a compact 
$S^1$-manifold~$M$, with action $S^1 \times M \to M$, $(s,x) \mapsto s \1 x$,
one proceeds as follows.
Given times $t_1 < \dots < t_{N} \in \RR$ and angles $s_1, \dots, s_N \in S^1$
one considers the functions
\begin{equation} \label{f:Chekanov}
f_t(x) \,=\,
\left\{ 
\begin{array}{ll}
f(x)     & \text{ if }\, t < t_1,  \\ [0.2em]
f(s_1 x) & \text{ if }\, t_1 \leqslant t < t_2, \\ [0.2em]
f((s_2 + s_1) \1 x) & \text{ if }\, t_2 \leqslant t < t_3, \\
\qquad \vdots & \\
f((s_N + \cdots +s_1) \1 x) & \text{ if }\, t_{N} \leqslant t,
\end{array}
\right.
\end{equation} 
and counts gradient ``$N$-jump flow lines'' of the vector field $- \nabla f_t$.
A neat way to see that a point $(t_1, \dots, t_{N}, s_1, \dots , s_N)$
corresponds to a point in $S^{2N-1}$ is through the join construction, 
\cite[\S 2.5]{BouOan12}.

This construction of equivariant Morse and Floer homology was explained in several
lectures by Chekanov~\cite{Che04}, and worked out by Noetzel~\cite{Noe08}, though never written up.
The construction and an isomorphism to the Borel construction is worked out for Morse homology in~\cite{Ber09}
based on~\cite{Che04},
and for Floer homology in~\cite{BouOan12} building on~\cite[\S~8b]{Se08}.
The construction and the isomorphism in~\cite{BouOan12} can be adapted to Rabinowitz--Floer homology,
yielding a homology $\RFH^\TT_{\mbox{\tiny jump}}(\Sigma,M)$ isomorphic to $\RFH^\TT (\Sigma,M)$.
The vanishing of $\RFH^\TT(\Sigma,M)$ for displaceable~$\Sigma$ then follows from the 
vanishing of $\RFH^\TT_{\mbox{\tiny jump}}(\Sigma,M)$, 
which in turn follows as for the non-equivariant $\RFH (\Sigma,M)$ by a leafwise intersection argument, 
because the chain groups of~$\RFH^\TT_{\mbox{\tiny jump}}(\Sigma,M)$ are exactly the chain groups of~$\RFH (\Sigma,M)$.



\newpage

\appendix

\section*{Appendix on transversality}

\subsection*{Preface}

Consider a smooth Morse function~$f$ on a closed Riemannian manifold~$(X,g)$.
In order to define the Morse homology of~$X$,
one needs that the spaces of gradient flow lines between critical points are manifolds.
This can be achieved by perturbing the Riemannian metric~$g$, see~\cite[\S \,2.3.2]{Schw93}:
The space of perturbations is the space of Riemannian metrics on~$X$.

In Floer theory on a symplectic manifold~$(M,\omega)$,
the role of~$f$ and~$X$ is played by the action functional~$\cA^H$ and the loop space
$\cL = C^\infty (S^1,M)$, 
and one works with $L^2$-Riemannian metrics on~$\cL$ given by $\omega$-compatible almost complex 
structures~$J$ on~$M$:
For $v \in \cL$ and vector fields $\xi_1, \xi_2$ along~$v$, 
$$
\langle \xi_1, \xi_2 \rangle_v \,:=\, \int_0^1 \omega (v(t)) \bigl( \xi_1(t), J(v(t)) \, \xi_2(t) \bigr) \, dt .
$$
Floer's equation (for $H=0$, say) is hence local in~$M$: 
\begin{equation}  \label{e:Floer.local}
\pp_s v (s,t) + J \bigl( v(s,t),t \bigr) \, \pp_t v(s,t) \,=\, 0 .
\end{equation}
Denote by $\widehat{\cM}_J(c_-,c_+)$ the space of solutions of~\eqref{e:Floer.local} asymptotic to 
critical points~$c_-, c_+$ of~$\cA^H$,
and by $\widehat{\cM}_J$ their union, with the $C_{\loc}^\infty$-topology.
The traditional way to achieve transversality for the spaces~$\widehat \cM_J(c_-,c_+)$ 
is by perturbing~$J$ in the space $\cJ_M$ of $\omega$-compatible almost complex structures on~$M$,
or in the slightly larger space $\cJ_{M \times S^1}$ of $1$-periodic loops in~$\cJ_M$.

These are very small perturbation spaces:
The almost complex structures in~$\cJ_M$ or~$\cJ_{M \times S^1}$ are only ``$M$-dependent'' or 
``$M \times S^1$-dependent''.
A much larger space of perturbations, 
which is the precise analogue of the space of Riemannian metrics playing the role 
of perturbations in Morse homology, 
is the space~$\cJ_\cL$ of $\omega$-compatible ``$\cL$-dependent'' almost complex structures:
An element $J \in \cJ_\cL$ associates with each $v \in \cL$ a loop $J_v (t)$ of 
$\omega$-compatible almost complex structures on~$T_{v(t)}M$.
The $L^2$-Riemannian metric on~$\cL$ is now
$$
\langle \xi_1, \xi_2 \rangle_v \,:=\, \int_0^1 \omega (v(t)) \bigl( \xi_1(t), J_v(t)\, \xi_2(t) \bigr) \, dt .
$$
Hence, Floer's equation for $J \in \cJ_\cL$ is non-local in~$M$:
\begin{equation*} 
\pp_s v (s,t) + J_{v(s)}(t) \,\pp_t v(s,t) \,=\, 0 .
\end{equation*}
%

The two main features one requires for the moduli spaces $\widehat \cM_J$ are 
that they are compact and that $\widehat \cM_J(c_-,c_+)$ are cut out transversally.
Both, compactness and transversality, which implies that $\widehat \cM_J(c_-,c_+)$ are manifolds,
are needed for defining the boundary operator~$d$ and proving $d^2=0$,
as well as for proving that the resulting homology does not depend on the choice of~$J$.

The advantage of working with a small space $\cP_{\sma}$ of perturbations, 
such as $\cJ_M$ or~$\cJ_{M \times S^1}$, 
is that if one can prove compactness of one moduli space~$\widehat \cM_J$,
then for the same reason one has compactness of~$\widehat \cM_{J'}$ for all $J' \in \cP_{\sma}$.
The disadvantage of these small perturbation spaces is that it is sometimes hard, or even impossible, 
to achieve transversality for a generic set of $J' \in \cP_{\sma}$.

Conversely, 
the advantage of working with a large space $\cP_{\lar}$ of perturbations, 
such as $\cJ_\cL$ or the perturbation space used for transversality in the theory of~M-polyfolds~\cite[\S \,5.3]{HWZ14}, 
is that transversality can be achieved easily, or at least in an easier way.
The disadvantage of these large spaces is that even if one has compactness of~$\widehat \cM_J$ for one~$J$,
compactness of~$\widehat \cM_{J'}$ for nearby $J' \in \cP_{\lar}$ can fail!

Assume now that one knows that $\widehat \cM_J$ is compact for one 
$J \in \cJ_M$.\footnote{
This is the case for Rabinowitz--Floer homology on an exact symplectic manifold 
and any $J \in \cJ_{\con}$ if one looks only at the part of~$\widehat \cM_J$ 
between critical points with action in a fixed interval~$[a,b]$, 
see\cite[Corollary~3.3]{CieFra09}, 
and similarly holds for these partial moduli spaces for 
$S^1$-equivariant Rabinowitz--Floer homology and any $J \in \cJ_{\con}$,
cf.~\S~\ref{ss:RFH} and \S~\ref{ss:eqRFH}.
}
The goal of this appendix is to explain that in this situation one can simultaneously achieve 
compactness and transversality for generic $J' \in \cP_{\lar}$ sufficiently close to~$J$,
not for the whole spaces $\widehat \cM_{J'}$, but for a compact part of~$\widehat \cM_{J'}$ near~$\widehat \cM_J$.

The natural framework for carrying out this argument
is the scale calculus, which forms a central building block of  M-polyfold theory, \cite{HWZ14}. 
In this setting, the (completion of the) loop space~$\cL$, the Hilbert-bundle over~$\cL$,
and also the space of perturbations~$\cP_{\lar}$,
are replaced by decreasing sequences of spaces and bundles,
and the gradient of the action functional and its Hessian at critical points  
become a scale vector field and a scale symmetric bilinear form.
Since the scale calculus has not yet become common knowledge,
and since setting up the scale structures in our situation would obscure the main idea,
we will stay in the framework of classical analysis and follow~\cite{RobSal95,RobSal01},
where ``finite step scale structures'' are already implicit.

To bring out the argument clearly, we first explain it in Appendix~\ref{s:app.Morse} in detail
in the framework of Morse homology on a finite-dimensional but possibly non-compact manifold~$X$.
The arguments are chosen such that they extend to the infinite-dimensional situation of Floer homology 
in~\S~\ref{ss:conley}, 
up to the compactness of the nearby parts of the spaces $\widehat \cM_{J'}$, 
that follows from a compactness result for non-local perturbations of the 
Cauchy--Riemann operator proven in~\S~\ref{ss:nonlocalcomp}.
Other parts of the construction of this Floer homology 
(such as the Fredholm theory and unique continuation) will be worked out in~\cite{FraSch17}.

\section{Transversality for Morse homology on finite-dimensional manifolds}
\label{s:app.Morse}

\subsection{The set-up} 
Let $X$ be a $C^\infty$-smooth finite-dimensional manifold.
We do not assume that $X$ is compact or complete,
but for convenience we assume that $X$ is connected.
Choose a Riemannian metric~$g$ on~$X$.
Let $f \colon X \to \RR$ be a $C^\infty$-smooth Morse function.
Denote by $\nabla f = \nabla_g f$ the gradient vector field of~$f$ with respect to~$g$.
We do not assume that this vector field is complete. 
However, we make a compactness assumption on the ``gradient flow lines'' of~$\nabla f$:
A gradient flow line $x \in C^\infty (\RR,X)$ is a solution of the ordinary differential equation
$$
\dot x(s) \,=\, -\nabla f(x(s)), \quad s \in \RR.
$$
%
For $a,b \in \RR$ let $\cG_a^b$ be the space of gradient flow lines~$x$ with $a \leqslant f(x(s)) \leqslant b$ for all $s \in \RR$. 
Endow $C^\infty (\RR, X)$ with the $C^\infty_{\loc}$-topology.
\begin{itemize}
\item[\bf (A)]
For all $a \leqslant b$ the space $\cG_a^b$ is compact in~$C^\infty (\RR, X)$.
\end{itemize}
Let $X_a^b = \{ x \in X \mid a \leqslant f(x) \leqslant b \}$ and $\Crit f_a^b = \Crit f \cap X_a^b$.
Assumption~(A) implies that $\cG_a^b$ is the disjoint union of the sets $\cG (x_-,x_+)$ of gradient flow lines
between points $x_-,x_+ \in \Crit f_a^b$, see \cite[Lemma~2.1]{CieFra11}.
Note that critical points of~$f$ are gradient flow lines.
Assumption (A) also implies that $\Crit f_a^b$ is a finite set.
The set $\cG_a^b$ is therefore the finite union of the sets $\cG (x_-,x_+)$ with $x_\pm \in \Crit f_a^b$.
 
The spaces $\cG (x_-,x_+)$ are not manifolds in general, 
but this holds true for a close-by Riemannian metric,~\cite[\S \,2.3.2]{Schw93}.
For any other Riemannian metric~$g'$ the negative gradient vector field $-\nabla_{g'}f$
is a pseudo-gradient vector field for~$f$. 
Anticipating the perturbation vector fields in the M-polyfold theory,
we work with pseudo-gradient vector fields $-\nabla f+v$ that agree with~$-\nabla f$ near $\Crit f_a^b$, 
see the next paragraph.
Note that such pseudo-gradient vector fields bijectively correspond to Riemannian metrics~$g'$ that agree with~$g$
near~$\Crit f_a^b$. 

The key argument for this approach to transversality is given 
in~\ref{ss:a.compact}, where we show that for every sufficiently small 
perturbation~$v$ one can select a compact part 
$\cG_a^b (v,\cN) = \bigcup_{x_-,x_+} \cG (x_-,x_+,v,\cN)$ 
of $\cG_a^b(v) = \bigcup_{x_-,x_+} \cG (x_-,x_+,v)$ near~$\cG_a^b$.
In~\ref{ss:a.lines} and~\ref{ss:a.bundle} we describe a Hilbert bundle $\cH \longleftarrow \cE$
over the Hilbert manifold~$\cH$ of lines from~$x_-$ to~$x_+$ near~$\cG_a^b$,
and show that the spaces $\cG (x_-,x_+,v,\cN)$ 
can be viewed as the zero-set of suitable sections of this bundle.
In~\ref{ss:a.smooth} we use (a minor variation of) the Sard--Smale theorem to show that
for a generic set of perturbations, $\cG (x_-,x_+,v,\cN)$ is smooth.

To make Appendix~\ref{s:app.Morse} useful for Appendix~\ref{s:transfloer},
we use the local compactness of~$X$ only when invoking the Arzel\`a--Ascoli theorem, a tool that is replaced
in~\S~\ref{ss:conley} by the compactness theorem in Floer homology proven in~\S~\ref{ss:nonlocalcomp}.

\subsection{Perturbations} \label{ss:a.per}
Fix $a<b$.
Let $c_1, \dots, c_N$ be the elements of $\Crit f_a^b$.
For each~$i$ choose an open neighbourhood~$U_i$ of~$c_i$ such that for $i \neq j$ 
the closures of $U_i$ and~$U_j$ are disjoint. 
Set $U = \bigcup U_i$.

The evaluation $\cG_a^b \to X$, $x \mapsto x(0)$, is continuous.
By assumption~(A) the image $K := \left\{ x(0) \mid x \in \cG_a^b \right\} \subset X$ is thus compact.
Since $\nabla f$ is a continuous vector field on~$X$,
we find an open and bounded neighbourhood~$\cN$ of~$K$ in~$X$
and $\delta >0$ such that $\| \nabla f (x) \| \geqslant \delta$ for all $x \in \cN \setminus U$. 
Choose $k \in \NN$ such that $k \geqslant 2$ and $k > \ind (c_i) -\ind (c_j)$ for all $c_i,c_j \in \Crit f_a^b$.
The vector space $\cV^k$ of $C^k$-vector fields on~$X$ that vanish on~$U$  
and have bounded derivatives up to order~$k$,
endowed with the $C^k$-norm, is a Banach space.
Its subset
$$
\cV^k_f \,=\, \bigl\{ v \in \cV^k \mid df(-\nabla f+v) < 0 
\,\mbox{ on\1 $\cN \setminus U$} \bigr\} 
$$
contains the open ball $B_\delta := \left\{ v \in \cV^k \mid \|v\|_{C^k} < \delta \right\}$,
since $df(-\nabla f+v) =  \langle -\nabla f +v, \nabla f \rangle$.

\begin{figure}[h]
 \begin{center}
  \psfrag{K}{$K$}
  \psfrag{N}{$\cN$}
  \psfrag{Ui}{$U_i$}
  \psfrag{Uj}{$U_j$}
  \leavevmode\epsfbox{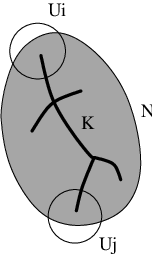}
 \end{center}
 \caption{The isolating neighbourhood $\cN$ of $K$.}
 \label{fig.NK}
\end{figure}
%
%

\subsection{Selection of compact perturbed moduli spaces} \label{ss:a.compact}
For $v \in \cV_f^k$ let $\cG_a^b(v,\overline \cN)$ be the set of solutions $x \colon \RR \to X$ of 
\begin{equation*}
\dot x (s) \,=\, -\nabla f(x(s)) + v (x(s)), 
\end{equation*}
such that $a \leqslant f(x(s)) \leqslant b$ and $x(s) \in \overline \cN$ for all $s \in \RR$.
The elements of $\cG_a^b(v,\overline \cN)$ are $C^{k+1}$-smooth, and the set of $j$th derivatives of elements of 
$\cG_a^b(v,\overline \cN)$ is uniformly bounded and equicontinuous for every $j \leqslant k$,
because $\overline \cN$ is bounded and because $\nabla f$ and $v$ are bounded on $\overline \cN$ 
with all their derivatives of order $\leqslant k$.
The Arzel\`a--Ascoli theorem thus implies that $\cG_a^b(v,\overline \cN)$ is $C^k_\loc$-compact.

The elements of $\cG_a^b(v,\overline \cN)$ may touch the boundary~$\pp \overline \cN$,
and even for ``regular'' $v$ (as defined in~\ref{ss:a.smooth}) 
this may lead to new ends of $\cG_a^b(v,\overline \cN)$
that do not come from breaking of flow lines and thus
may obstruct the property~$d^2=0$ necessary to define Morse homology.
As the next lemma shows, this cannot happen if $v$ is small enough.
Recall that $B_{\delta} \subset \cV_f^k$.

\begin{lemma} \label{le:epsf}
There exists $\eps \in (0,\delta)$ such that for all $v \in B_\eps$,
$$
x(\RR) \subset \cN \quad \mbox{ for all }\, x \in \cG_a^b(v,\overline \cN) .
$$
\end{lemma}

In other words, $\overline \cN$ is an isolating neighbourhood for the flow of $-\nabla f+v$
simultaneously for all~$v \in B_\eps$,
and for each such~$v$ the set $\cG_a^b(v,\overline \cN)$ is an isolated invariant set of this flow.

\proof
If not, there exists a sequence of perturbations $v_\nu \in \cV_f^k$ with $\lim_{\nu \to \infty} v_\nu \to 0$ in~$C^k$
and a sequence of flow lines $x_\nu \in \cG_a^b(v_\nu, \overline \cN)$
such that $x_\nu(s_\nu) \in \pp \cN$ for some time $s_\nu \in \RR$.
Since $\cG_a^b(v_\nu, \overline \cN)$ is invariant under time-shift, we can assume that
$s_\nu =0$ for all~$\nu$.
Recall that $x_\nu$ solves
$$
\dot x_\nu (s) \,=\, -\nabla f(x_\nu(s)) + v_\nu (x_\nu(s)), \quad s \in \RR .
$$
Since $x_\nu (\RR) \subset \overline \cN$ and since the sequence $v_\nu$ is uniformly bounded in~$C^k$,
we can invoke again the Arzel\`a--Ascoli theorem to find a subsequence that converges in
$C^k_\loc$ to an element $x \in \cG_a^b$. Hence $x(0) \in K$.
On the other hand, $x(0) = \lim x_{\nu_j} (0) \in \pp \cN$.
But $K \cap \pp \cN = \emptyset$, a contradiction.
\proofend

\subsection{The Hilbert manifold of lines} \label{ss:a.lines}
Fix two points $x_-, x_+ \in \cN$.
Consider the ``set of lines'' from $x_-$ to~$x_+$ in~$\cN$,
$$
\cH \,:=\, \cH_{x_-,x_+}(\cN) \,:=\, \mbox{
``$\bigl\{ x \in W^{1,2}(\RR,\cN) \mid \displaystyle \lim_{s \to \pm \infty} x(s) = x_\pm \bigr\}$''} .
$$
We describe this set by constructing a Hilbert manifold atlas: 
Let $C^\infty_c(\RR,\cN,x_-,x_+)$ be the set of $x \in C^\infty (\RR,\cN)$ such that there exists $T>0$ with 
$$
x(s) = x_- \; \mbox { for } s \leqslant -T, \qquad x(s) = x_+ \; \mbox { for } s \geqslant T.
$$
Fix $x \in C^\infty_c(\RR,\cN,x_-,x_+)$.
The pull-back bundle $x^* T\cN \to \RR$ is a trivializable vector bundle over~$\RR$ of dimension $n = \dim \cN$. 
Fix a trivialization $\Phi \colon x^* T\cN \to \RR \times \RR^n$.
Choose an open neighbourhood $V_x \subset x^*T\cN$ of the zero-section such that
the exponential map $\exp_g \colon T_{x(r)}\cN \cap V_x \to \cN$ is injective for all $r \in \RR$.
The set 
$$
U_x \,:=\, 
\left\{ \xi \in W^{1,2} (\RR,\RR^n) \mid (r,\xi(r)) \in \Phi (V_x) \mbox{ for all } r \in \RR \right\}
$$
is an open subset of $W^{1,2}(\RR,\RR^n)$.
Define the injective map
$$
\phi_x \colon U_x \to C^0(\RR,\cN,x_-,x_+), \quad \xi \mapsto \exp_g \left( \Phi^{-1}(\xi) \right) ,
$$
where $C^0(\RR,\cN,x_-,x_+)$ is the set of continuous maps $y \colon \RR \to \cN$ with $\lim_{s \to \pm \infty} y(s) = x_\pm$.  
Using elementary Sobolev theory one shows that the transition maps 
$(\phi_{x'})^{-1} \circ \phi_x \colon U_{x'} \cap U_x \to U_{x'} \cap U_x$ are smooth in~$W^{1,2}(\RR, \RR^n)$.
Now define 
$$
\cH \,:=\, \bigcup_{x \in C^\infty_c(\RR,\cN,x_-,x_+)} \phi_x (U_x) .
$$
The space $\cH$ is thus a Hilbert manifold modeled on~$W^{1,2}(\RR, \RR^n)$.

\subsection{The Hilbert bundle set-up}  \label{ss:a.bundle}
The tangent space at $x \in \cH$ is
$$
T_x \cH \,=\, W^{1,2}(\RR, x^* T\cN) .
$$
A larger bundle over~$\cH$ is the $L^2$-bundle with fibre $\cE_x = L^2 (\RR, x^*T\cN)$,

$$
\xymatrix{ 
T\cH \; \ar@{^{(}->}[r] \ar[d] &
\, \cE \ar[ld] & 
\\
\cH \,& 
}
$$

%
%
%

\ni
Fix now $x_-,x_+ \in \Crit f$.
The vector fields $-\nabla f+v$ with $v \in \cV^k_f$ are $C^k$-smooth pseudo-gradient vector fields for~$f$ 
that agree with $-\nabla f$ on~$U \supset \Crit f_a^b$. 
For $v \in \cV_f^k$ denote by $\cG(x_-,x_+,v,\cN)$ the set of solutions $x \colon \RR \to \cN$ of 
\begin{equation} \label{e:gradh}
\dot x(s) \,=\, (-\nabla f +v) (x(s)) .
\end{equation}
We can now interpret the set $\cG (x_-,x_+,v,\cN)$ as the zero-locus of a section:
Recall that a section of the bundle $\cH \stackrel{p}{\longleftarrow} \cE$ is a map 
$\cH \stackrel{s}{\longrightarrow} \cE$ such that $p \circ s = \id_\cH$.
To $v \in \cV_f^k$ associate the $C^k$-section $s_v \colon \cH \to  \cE$ defined by
$$
s_v(x) \,:=\, \dot x + (\nabla f -v) (x) .
$$
Then $x \in s_v^{-1}(0)$ if and only if $x$ is a $W^{1,2}$-flow line of $-\nabla f +v$ in~$\cN$ from $x_-$ to~$x_+$.
Using equation~\eqref{e:gradh} we see that such an~$x$ is actually $C^{k+1}$-smooth.

For $x$ in the zero-section~$\cH$ of the bundle $\cE \stackrel{p}{\longrightarrow} \cH$ we have 
the canonical splitting
$$
T_x \cE \,=\, \cE_x \oplus T_x \cH .
$$ 
(Defining such a splitting at a point~$x$ off the zero-section would require the choice of a connection.)
The differential $d s_v (x)$ is a map $T_x \cH \to T_{s_v(x)} \cE$.
For $x \in s_v^{-1}(0)$ define the {\it vertical differential}\/ by
$$
D s_v(x) \colon T_x \cH \to \cE_x,  \quad Ds_v(x) = \pi \circ ds_v (x) 
$$
where $\pi \colon T_x \cE = \cE_x \oplus T_x \cH \to \cE_x$ is the projection along~$T_x \cH$.

A $C^1$-section $s \colon \cH \to \cE$ is said to be {\it transverse to the zero-section}, $s \pitchfork 0$,
if $Ds(x)$ is surjective for all $x \in s^{-1}(0)$.
Furthermore, a $C^1$-section $s \colon \cH \to \cE$ is called {\it Fredholm}\/
if $Ds(x) \colon T_x \cH \to \cE_x$ is a Fredholm operator for all $x \in s^{-1}(0)$.

\begin{proposition} \label{p:Satz1}
Let $s \colon \cH \to \cE$ be a $C^k$-smooth Fredholm section such that $s \pitchfork 0$.
Then $s^{-1}(0)$ is a $C^k$-smooth manifold, and for $x \in s^{-1}(0)$ we have
$$
\dim_x s^{-1}(0) \,=\,  \dim \ker Ds(x) \,=\, \ind \, Ds(x).
$$
Here $\dim_x s^{-1}(0)$ denotes the dimension of the connected component containing~$x$, 
and $\ind \, Ds(x) := \dim \ker Ds(x) - \dim \coker \, Ds(x)$ is the Fredholm index of~$Ds(x)$.
\end{proposition} 

The first identity follows from the implicit function theorem for $C^k$-maps between Banach manifolds,
and the second identity holds because $Ds(x)$ is surjective, i.e., $\coker \, Ds(x) = 0$.

Recall that the support of every perturbation $v \in \cV_f^k$ of $-\nabla f$ is disjoint from $\Crit f_a^b$.
The sections~$s_v$ are thus Fredholm, see \cite[Theorem~2.1]{RobSal95} or~\cite[\S \,2.2]{Schw93}.
These references also prove

\begin{corollary}  \label{cor:index}
If $x_-, x_+$ have Morse index \1$\ind (x_\pm)$
and if $s_v \pitchfork 0$, then  
$$
\dim_x s_v^{-1}(0) \,=\, \ind (x_-) - \ind (x_+) .
$$
\end{corollary}

\subsection{Finding smooth perturbed moduli spaces} \label{ss:a.smooth}
Let $\eps \in (0,\delta)$ be as in Lemma~\ref{le:epsf}.
The open subset~$B_\eps$ of the Banach space~$\cV^k$ is homeomorphic to~$\cV^k$ and hence metrizable by a complete metric. 
In particular, $B_\eps$ is a Baire space, i.e., every subset of~$B_\eps$ of the second category 
is dense in~$B_\eps$.
%

Fix $x_-,x_+ \in \Crit f_a^b$ and write again $\cH = \cH_{x_-,x_+}(\cN)$.
Consider the $C^k$-map 
\begin{equation} \label{e:S}
S \colon B_\eps \times \cH \to \cE, \quad S(v,x) = s_v(x) .
\end{equation}
By Lemma~\ref{le:epsf} the set $\cG (x_-,x_+,v,\cN)$ 
is the set of flow lines in $\cG_a^b(v,\overline \cN)$ from $x_-$ to~$x_+$. 
Write 
$$
\cG(x_-,x_+,\cN) \,=\, \bigcup_{v \in B_\eps} \{v\} \times \cG(x_-,x_+,v,\cN) .
$$
Then $\cG (x_-, x_+,\cN) = S^{-1} (0)$.  

\begin{theorem} \label{t:trans}
There exists a subset $\cV_\reg (x_-,x_+) \subset B_\eps$ of the second category such that 
$s_v \colon \{v\} \times \cH \to \cE$ is tranverse to~$0$ for all $v \in \cV_\reg (x_-,x_+)$.
\end{theorem}

Before giving the proof, we draw the desired conclusions.
Recall that the set $B_\eps$ depends only on~$f$ and $a,b$, 
while the set $\cV_\reg (x_-,x_+) \subset B_\eps$ in Theorem~\ref{t:trans} depends also on $x_-,x_+$.
Proposition~\ref{p:Satz1}, Corollary~\ref{cor:index} and Theorem~\ref{t:trans} show that for all~$v$ 
in the subset $\cV_\reg (x_-,x_+) \subset B_\eps$ of the second category, 
the set $s_v^{-1}(0) = \{v\} \times \cG(x_-,x_+,v,\cN)$ is a $C^k$-smooth manifold of dimension 
$\ind (x_-) - \ind (x_+)$.
Recall that there are only finitely many critical points in $\Crit f_a^b$.
With $\cV_\reg := \bigcap_{x_-,x_+} \cV_\reg (x_-,x_+)$ we find

\begin{corollary} \label{cor:trans}
Assume that (A) holds.
Then there exists a subset $\cV_{\reg} \subset B_\eps$ of the second category such that 
for every $v \in \cV_{\reg}$ all the sets 
$\cG (x_-,x_+,v,\cN)$ with $x_\pm \in \Crit f_a^b$
are $C^k$-manifolds of dimension $\ind (x_-) - \ind (x_+)$. 
\end{corollary}

\begin{remark}
{\rm
An argument due to Taubes, for which we refer to~\cite[\S \, 3.2]{McSal12}, now implies that Corollary~\ref{cor:trans}
holds with $k$ replaces by~$\infty$.
We do not need this in the sequel.
}
\end{remark}

Recall that the sets $\cG_a^b(v,\overline \cN)$ are $C_\loc^k$-compact for all $v \in \cV_f^k$.
As in the proof of Lemma~2.1 in~\cite{CieFra11} if follows that for these~$v$,
%
\begin{equation} \label{e:coprod}
\cG_a^b(v,\overline \cN) \,=\, \coprod_{x_-,x_+} \cG (x_-,x_+,v, \overline \cN),
\qquad  x_\pm \in \Crit f_a^b .
\end{equation}
By Lemma~\ref{le:epsf} we have $\cG (x_-,x_+,v, \overline \cN) = \cG (x_-,x_+,v,\cN)$ for $v \in B_\eps$,
and
by Corollary~\ref{cor:trans} these sets are $C^k$-manifolds 
in~$\cN$ of dimension $\ind (x_-) - \ind (x_+)$ for every $v \in \cV_\reg$.
Morse homology can now be defined by using the spaces~\eqref{e:coprod} with $v \in \cV_\reg$.

\m \ni
{\it Proof of Theorem~\ref{t:trans}.}
We invoke an abstract result, that is a minor variation
of the classical Sard--Smale theorem, proved by Sard and Quinn~\cite[Theorem~1]{QuSa72}
and put in the form below by Henry~\cite[Theorem~5.4]{Hen05}, see also~\cite{Sch13}.
In the formulation we use the letters from the previous text, but change the fonts of the spaces 
to stress the abstract setting.
A map $f \colon X \to Y$ between topological spaces is called $\sigma$-proper if $X$ can be written as a countable union
$X = \bigcup_i X_i$ such that the restrictions $f \colon X_i \to Y$ are proper.

\begin{proposition} \label{p:abstract}
Let $\mathcal V$ and $\mathcal H$ be Banach manifolds
and let $\mathcal{E} \to \mathcal{H}$ be a Banach bundle.
Let $S \colon \mathcal V \times \mathcal H \to \mathcal E$ be a $C^k$-map
such that $s_v = S(v, \cdot) \colon \{v\} \times {\mathcal H} \to \mathcal E$
is a section for each $v \in \mathcal{V}$. 
Assume that
\begin{itemize}
\item[(i)]
For every $(v,x) \in S^{-1}(0)$ the vertical differential $DS (v,x)$ is surjective.

\m
\item[(ii)]
For every $v \in \mathcal V$ the section $s_v \colon \{v\} \times {\mathcal H} \to \mathcal{E}$
is a Fredholm map of index $<k$.

\m
\item[(iii)]
The projection $S^{-1}(0) \to \mathcal{V}$, $(v,x) \mapsto v$, is $\sigma$-proper.
\end{itemize}
Then there exists a subset $\mathcal{V}_\reg \subset \mathcal{V}$ of the second category such that 
$s_v \pitchfork 0$ for all $v \in \mathcal{V}_\reg$.
\end{proposition}

In this version, the separability assumption on $\mathcal{V}$ and~$\mathcal{H}$ in the classical Sard--Smale theo\-rem
is replaced by the $\sigma$-properness assumption~(iii).
While this modification is irrelevant in the case of finite-dimensional Morse homology, 
in which $B_\eps$ and $\cH$ are separable, 
it will be relevant in the infinite-dimensional situation of Floer homology addressed in~\S~\ref{ss:conley},
where neither $\mathcal{V}$ nor~$\mathcal{H}$ is separable anymore.

We apply Proposition~\ref{p:abstract} with $\mathcal{V} = B_\eps$,
$\mathcal{H} = \cH$, $\mathcal{E} = \cE$, 
and $S$ as in~\eqref{e:S}, so that $S^{-1}(0) = \cG(x_-,x_+, \cN)$.

We first verify assumption~(iii) of Proposition~\ref{p:abstract}.
Let $U$ be the neighbourhood of~$\Crit f_a^b$ defined in~\S~\ref{ss:a.per}.
After choosing $U$ smaller, if necessary, we can assume that the Morse Lemma
holds for~$f$ on~$\overline U$.
For $T \in \NN$ define 
$$
\cG_T(x_-,x_+, \cN) \,=\, 
\left\{ (v,x) \in S^{-1}(0) \mid x(s) \in \overline U \1\mbox{ for }\1 s \notin [-T,T] \right\}.
$$
In other words, $(v,x) \in S^{-1}(0) = \cG(x_-,x_+,\cN)$  belongs to~$\cG_T(x_-,x_+,\cN)$ 
if $x(s) \in \overline U(x_-)$ for $s < -T$ and $x(s) \in \overline U(x_+)$ for $s > T$.
Then $S^{-1}(0) = \bigcup_{T \in \NN} \cG_T(x_-,x_+,\cN)$.
We thus need to prove

\begin{lemma}  \label{le:compact}
For every $T \in \NN$ the projection $p \colon \cG_T(x_-,x_+,\cN) \to B_\eps$, $(v,x) \mapsto v$, is proper.
\end{lemma}
 
\proof
The reason for this is that the only possible source of non-compactness of a sequence 
$(v_\nu, x_\nu) \subset \cG_T(x_-,x_+,\cN)$ with $(v_\nu)$ in a compact set of~$B_\eps$ 
is breaking, which does not occur because $x_\nu(s) \subset \overline U(x_-)$ for $s < -T$ 
and $x_\nu(s) \subset \overline U(x_+)$ for $s > T$ and because $f$ is Morse.

Given a compact subset $K \subset B_\eps$ we wish to show that $p^{-1}(K) \subset \cG_T(x_-,x_+,\cN)$ is compact. 
Let $(v_\nu, x_\nu)$ be a sequence in~$p^{-1}(K)$.
Since $K$ is compact and by the Arzel\`a--Ascoli theorem we can pass to a subsequence such that
$v_\nu \to v$ in~$B_\eps$ and $x_\nu|_{[-T,T]} \to x_{[-T,T]}$ in~$C^k([-T,T],\overline \cN)$.
Recall that on~$\overline U$ we have $v_\nu = v = 0$, and that the Morse Lemma holds for~$f$ on~$\overline U$. 
Since $x_\nu(\pm T) \to x_{[-T,T]}(\pm T)$, we also have 
$x_\nu|_{(-\infty,-T]} \to x_{(-\infty,-T]}$ in~$C^k((-\infty,-T], \overline \cN)$ and
$x_\nu|_{[T,\infty)} \to x_{[T,\infty)}$ in~$C^k([T,\infty), \overline \cN)$, 
with exponential decay of all derivatives of order $\leqslant k$, by the Morse Lemma. 
For the concatenation~$x = x_{(-\infty,-T]} \# x_{[-T,T]} \# x_{[T,\infty)}$
we thus have $x_\nu \to x$ in $C^k(\RR,\overline \cN)$,
with exponential decay of all derivatives of order $\leqslant k$.
By the Arzel\`a--Ascoli theorem and by Lemma~\ref{le:epsf}, $(v,x) \in \cG_T(x_-,x_+,\cN)$.
Furthermore, the convergence $x_\nu \to x$ in $C^1(\RR,\overline \cN)$ with exponential decay 
of all derivatives of order $\leqslant 1$
implies that $(v_\nu,x_\nu) \to (v,x)$ in $B_\eps \times \cH$.
\proofend

That assertions~(i) and (ii) of Proposition~\ref{p:abstract} hold in our setting is a standard result
(that in fact holds for the unrestricted map $S \colon \cV_f^k \times \cH \to \cE$).
We have already seen this for assertion~(ii).
It thus remains to prove

\begin{proposition} \label{p:DSsur}
For every $(v,x) \in S^{-1}(0)$ the vertical differential $DS (v,x)$ is surjective.
\end{proposition}

\proof
We distinguish two cases.

\smallskip \noindent
{\bf Case 1. $x_- \neq x_+$.}
Fix $(v,x) \in S^{-1}(0)$. Then $\dot x + (\nabla f -v) (x) =0$.
The tangent space $T_v B_\eps$ is the Banach space~$\cV^k$.  
The vertical differential $DS(v,x) \colon \cV^k \oplus T_x \cH \to \cE_x$ is given by
\begin{equation} \label{e:DS}
DS(v,x)(\hat v, \hat x) \,=\, Ds_v(x) \hat x + \hat v .
\end{equation}
In particular, the image of $Ds_v(x)$ is contained in the image of $DS(v,x)$.
Since $Ds_v(x)$ is Fredholm, its image is closed in~$\cE_x$ and of finite codimension, 
and hence the image of~$DS(v,x)$ is also closed in~$\cE_x$.

The fiber~$\cE_x$ is equipped with the $L^2$-inner product  
$$
\langle \eta_1, \eta_2 \rangle_g \,:=\, \int_{-\infty}^\infty g_{x(s)} \left( \eta_1(s), \eta_2(s) \right) ds.
$$
Assume that $\eta \in \coker \, DS(v,x) \subset \cE_x$, that is, 
$\eta$ is $\langle \, , \, \rangle_g$-orthogonal to the image of~$DS(v,x)$.
We must show that $\eta = 0$.
By assumption, 
$$
\langle DS (v,x) (\hat v, \hat x), \eta \rangle_g =0 
\quad \mbox{ for all }\, (\hat v, \hat x) \in \cV^k \oplus T_x \cH .
$$
In particular
\begin{eqnarray}
\langle Ds_v(x) \hat x, \eta \rangle_g \,=\, 0 && \mbox{ for all }\, \hat x \in T_x \cH ,   \label{e:ds1}  \\ 
\langle \hat v, \eta \rangle_g \,=\, 0         && \mbox{ for all }\, \hat v \in \cV^k . \label{e:ds2} 
\end{eqnarray}
By~\eqref{e:ds1}, $\eta \in \coker \, Ds_v(x) = \ker Ds_v(x)^*$,
where $Ds_v(x)^*$ denotes the adjoint operator. 
Choose a chart $\RR^n \to \cN$ covering the whole flow line~$x$.
In this chart, 
$$
D s_v (x) \,\hat x \,=\, \dot {\hat x} + \Hess f(x) \,\hat x - dv (x) \,\hat x.
$$ 
Integrating by parts, we find that the adjoint operator is given by
$$
Ds_v(x)^* \,\eta \,=\, - \dot \eta + \bigl[ \Hess f(x) - dv (x) \bigr]^T \eta .
$$
An element $\eta \in \coker \, Ds_v(x) = \ker Ds_v(x)^*$ is thus a solution of the ordinary differential equation
with $C^{k-1}$-coefficients
\begin{equation} \label{e:ODECk-1}
\dot \eta (s) \,=\, \bigl[ \Hess f(x(s)) - dh^T (x(s)) \bigr] \, \eta (s), \quad s \in \RR .
\end{equation}
Hence $\eta \in C^k(\RR,x^*T \cN)$.
We must show that $\eta (s) = 0$ for all $s \in \RR$.

Recall that the sets $\overline U_i$ are mutually disjoint.
Since $x_- \neq x_+$, the flow line~$x$ therefore intersects the set $W := \cN \setminus \overline U$.
Assume that there exists $s^* \in \RR$ such that $x^* := x(s^*) \in W$ and $\eta(s^*) \neq 0$.
Since $\eta$ is continuous, there exists $\eps >0$ such that
$\eta (s) \neq 0$ for all $s \in [s^* - \eps, s^*+\eps]$.
Since $v \in B_\eps \subset \cV_f^k$, $x$ is a non-constant flow line of $-\nabla f +v$,
and therefore the function $f \circ x$ is strictly decreasing. 
Hence the set
$$
W_\eps \,:=\, \left\{ x \in \cN \mid f(x(s^*+\eps)) < f(x) < f(x(s^*-\eps)) \right\} \cap W
$$
is an open neighbourhood of~$x^*$, and
\begin{equation} \label{e:We}
x(s) \notin W_\eps \quad \mbox{ for $s \notin (s^*-\eps , s^*+\eps)$.}
\end{equation}
Choose $\hat v \in \cV^k$ with support in $W_\eps$ such that
$\langle \hat v (x^*), \eta (x^*) \rangle >0$ and $\langle \hat v (x(s)), \eta (x(s)) \rangle \geqslant 0$ 
for all $s \in (s^*-\eps , s^*-\eps)$.
By~\eqref{e:We}, $\langle \hat v (x(s)), \eta (x(s)) \rangle =0$ for $s \notin (s^*-\eps , s^*-\eps)$.
Therefore, 
$$
\langle \hat v, \eta \rangle_g \,:=\, 
\int_{-\infty}^\infty \big\langle \hat v(x(s)), \eta (s) \big\rangle \,ds \,=\,
\int_{s^*-\eps}^{s^*+\eps} \big\langle \hat v(x(s)), \eta (s) \big\rangle \,ds 
\,>\, 0,
$$
in contradiction to~\eqref{e:ds2}.

We are left with showing that $\eta (s) = 0$ also for $x(s) \in \overline U$.
Assume that $x(s^*) \in \overline U_i$.
Since the flow line~$x$ is not entirely contained in~$\overline U_i$,
there exists $s'$ with $x(s') \in W$. We have already seen that $\eta (s) = 0$ for $s$ near~$s'$.
Since $\eta$ is a solution of the ordinary differential equation~\eqref{e:ODECk-1} with $C^{k-1}$-coefficients,
and since $k \geqslant 2$,
it follows that $\eta (s') =0$.
We have shown that $\eta (s)=0$ for all $s \in \RR$.

\medskip \noindent
{\bf Case 2. $x_- = x_+$.}
In this case, $S(v,x) = s_v(x) = \dot x + \nabla f(x) =0$ if and only if $x \equiv x_0 := x_-=x_+$
is the constant flow line.
Recall from Case~1 that the image of~$Ds_v(x)$ is contained in the image of~$DS(v,x)$.
It thus suffices to show that for every $v \in B_\eps$ the operator $Ds_v(x_0) \colon T_{x_0}\cH \to \cE_{x_0}$
is surjective, or, equivalently, that the adjoint operator $Ds_v(x_0)^*$ is injective.
As in~\eqref{e:ODECk-1} an element $\eta \in \ker Ds_v (x_0)^*$ is a solution of
\begin{equation} \label{e:DA}
\dot \eta (s) \,=\, A \, \eta (s), \quad s \in \RR,  
\end{equation}
where $A = \Hess f (x_0)$.
In a suitable basis, $A = \diag (a_1, \dots, a_n)$ with $a_i \neq 0$ because $f$ is Morse. 
The solutions of~\eqref{e:DA} are of the form $\eta (s) = \left( \eta_1(s), \dots, \eta_n(s) \right)$
with $\eta_i(s) = \eta_i(0) e^{a_i s}$. They lie in $L^2$ only if $\eta_i(0) =0$, i.e., if $\eta =0$.
\proofend

\section{Transversality for Floer homology with one compact moduli space} 
\label{s:transfloer}

\subsection{The Conley-type argument in the Floer case} \label{ss:conley}

In this paragraph we outline how the transversality scheme of Appendix~\ref{s:app.Morse},
and in particular the Conley-type argument in Section~\ref{ss:a.compact},
works in the setting of a Floer homology for which one has compactness of the moduli space~$\widehat \cM_J$ 
for one almost complex structure~$J$.
This is, for instance, the case for exact Liouville domains
or for symplectically aspherical closed symplectic manifolds.
The key difference to the Morse homology on a finite-dimensional manifold~$X$
is that now the role of~$X$ is taken by an infinite-dimensional manifold,
that is not locally compact.
The Arzel\`a--Ascoli theorem is thus not available anymore;
it is substituted by a compactness theorem for solutions of the Cauchy--Riemann equation 
with a non-local perturbation, that we address in the next paragraph.

\m
\ni 
{\bf Dictionary.}
Let $(M,\omega)$ be a symplectic manifold, choose an $S^1$-family~$J_t$ of almost complex structures
compatible with~$\omega$, and 
let $H \colon M \times S^1 \to \RR$ be a smooth function such that all 1-periodic
orbits of its Hamiltonian field~$X_{H_t}$ are non-degenerate.
To avoid technical and notational complications, we assume that $M=\CC^n$.
For $\ell \geqslant 0$ consider the Sobolev-spaces of loops $H_\ell = W^{\ell,2}(S^1, \CC^n)$. 
By assumption~(A) below the set 
$\Crit \cA_a^b = \{c_1, \dots, c_N\}$ of critical points of the 
action functional $\cA := \cA_H$ with action in $[a,b]$ is finite. 
Choose $k \in \NN$ such that $k \geqslant 2$ and $k > \ind (c_i) -\ind (c_j)$ for all $c_i,c_j \in \Crit \cA_a^b$.

The role of the manifold~$X$ and of the Morse function $f \colon X \to \RR$ is now taken by~$H_k$ and 
the action functional $\cA \colon H_k \to \RR$, 
and the gradient flow line equation $\dot x (s) = -\nabla f (x(s))$ becomes $\dot x (s) = -\nabla \cA (x(s))$. 
If one takes $\nabla$ to be the $H_0$-gradient with respect to the 
Riemannian metric on~$H_k$ induced by $\omega$ and~$J_t$, 
this is Floer's equation:
A gradient flow line $x \colon \RR \to H_k$ is a solution of the ordinary differential equation on~$H_k$
\begin{equation*} 
\dot x(s) \,=\, 
- J_t(x(s)(t))  \bigl( \partial_t x(s) + X_{H_t}( x(s)(t) ) \bigr) , 
\quad s \in \RR.
\end{equation*}
In the sequel, we omit the lower order term $-J_t \,X_{H_t}$ from the notation,
and for the situation of ($S^1$-equivariant) Rabinowitz--Floer homology we also neglect the $\RR$-factor of~$\eta$
and the spheres~$S^{2N+1}$.
For $a,b \in \RR$ let $\mathcal{G}_a^b$ be the space of gradient flow lines~$x$ with $a \leqslant \cA(x(s)) \leqslant b$ for all $s \in \RR$. 
Endow $C^\infty (\RR, H_k)$ with the $C^\infty_{\loc}$-topology.
We again assume that
\begin{itemize}
\item[\bf (A)]
For all $a \leqslant b$ the space $\mathcal{G}_a^b$ is compact in~$C^\infty (\RR, H_k)$.
\end{itemize}
This assumption is satisfied if $M$ is compact and $[\omega]$ vanishes on~$\pi_2(M)$, 
such as exact Liouville domains, 
and also for the gradient flow lines with cascades used to define ($S^1$-equivariant) Rabinowitz--Floer homology.

Choose open neighbourhoods~$U_i$ of~$c_i$ in~$H_1$ with disjoint closures.
The evaluation $\mathcal{G}_a^b \to H_k$ at~$0$ is continuous and has compact image~$K \subset H_k$, 
and $\nabla \cA \colon H_k \to H_{k-1}$ is smooth and the inclusion $H_{k-1} \to H_0$ is continuous. 
We can thus find an open and bounded neighbourhood~$\mathcal N$ of~$K$ in~$H_k$ and 
$\delta >0$ such that
$\| \nabla \cA \|_{H_0} \geqslant \delta$ on~$\mathcal N \setminus U$.

Recall that given two Banach spaces $(Y, \|\;\|_Y)$, $(Z,\|\;\|_Z)$
that are continuously embedded in the same Banach space, 
the space $Y \cap Z$ with norm $\|\;\|_Y + \|\;\|_Z$ is a Banach space. 
For $j \geqslant 0$ let $C^k(H_j,H_j)$ be the (non-separable) Banach space of $C^k$-vector fields $H_j \to H_j$ 
that have bounded derivatives up to order~$k$, endowed with the $C^k$-norm.
For every $j \in \{ 1, \dots, k \}$ the restriction from $H_j$ to~$H_k$ 
defines a continuous map $C^k(H_j,H_j) \to C^0 (H_k,H_1)$,
which is an embedding since the inclusion $H_k \subset H_j$ is dense.
With respect to these embeddings define the Banach space
\begin{equation} \label{e:CkH}
\mathcal{C}^k(H,H) \,=\, C^k(H_k,H_k) \cap C^k(H_{k-1},H_{k-1}) \cap \dots \cap C^k(H_1,H_1) .
\end{equation}
For instance, $\mathcal{C}^2(H,H) =  C^2(H_2,H_2) \cap C^2(H_1,H_1)$ consists of all those $C^2$-maps $H_1 \to H_1$ 
whose restriction to~$H_2$ takes values in $H_2$, and is still $C^2$ as a map $H_2 \to H_2$.

Now define $\mathcal{V}^k$ to be the closed subspace of $\mathcal{C}^k(H,H)$
formed by those vector fields that vanish on~$U$.
The elements $\hat v \in \mathcal{V}^k$ used to extend the proof of Proposition~\ref{p:DSsur}
are ``bump vector fields'', that in a chart around a point in~$H_1  \setminus \overline U$
are of the form
$$
x \mapsto \beta^\eps(\|x\|_1) \,e_i
$$
where $\eps >0$ and $\beta^\eps \colon \RR \to [0,\eps]$ is a smooth function that has support in $[-\eps, \eps]$
and that is equal to~$\eps$ near~$0$, 
and where $e_i$ is a vector of the Fourier basis of $H_1$.
The subset
$$
\mathcal{V}^k_{\cA} \,=\, 
\bigl\{ v \in \mathcal{V}^k \mid d\cA \2 (-\nabla \cA + v) < 0 
\,\mbox{ on\1 $\mathcal N \setminus U$} \bigr\} 
$$
of $\mathcal{V}^k$ contains the open $\delta$-ball in~$\mathcal V^k$,
since $d\cA \2 (-\nabla \cA + v) = \langle -\nabla \cA +v, \nabla \cA \rangle_{H_0}$.
The Hilbert manifold $\mathcal{H}_{x_-,x_+}(\mathcal N)$ is constructed as in~\ref{ss:a.lines},
but is now modeled on $W^{1,2}(\RR, H_k)$ and is therefore not separable.
The Hilbert bundle $\mathcal{E} \to \mathcal{H}$ has fibre $\mathcal{E}_x = L^2(\RR,H_{k-1})$.

For $v \in \mathcal{V}_{\cA}^k$ let $\mathcal{G}_a^b (v,\overline{\mathcal{N}})$ 
be the set of solutions $x \colon \RR \to H_k$ of 
\begin{equation*}
\dot x (s) \,=\, -\nabla \cA (x(s)) + v (x(s))
\end{equation*}
such that $a \leqslant \cA (x(s)) \leqslant b$ and $x(s) \in \overline{\mathcal{N}}$ for all $s \in \RR$.
The main difference to the case of Morse homology on a finite-dimensional manifold studied in Appendix~\ref{s:app.Morse}
is that now $H_k$ is not locally compact,
whence we cannot appeal to the Arzel\`a--Ascoli theorem anymore.
It is now Theorem~\ref{t:compact} below that shows that for $v \in B_\delta$ the spaces 
$\mathcal{G}_a^b(v, \overline{\mathcal N})$ are $C_{\loc}^0$-compact in~$C^0(\RR,H_k)$. 
Since $v \in \cV^k$ we can then use bootstrapping to see that 
$\mathcal{G}_a^b(v,\mathcal N)$ are $C_{\loc}^k$-compact in~$C^k(\RR,H_k)$.
We can now follow Appendix~\ref{s:app.Morse} to obtain transversality 
for the perturbed moduli spaces~$\mathcal{G}_a^b(v,\mathcal N)$
for a generic set of perturbations~$\mathcal{V}_{\reg} \subset B_\eps \subset B_\delta$. 
Some of the further tools needed are:
\begin{itemize}
\item[$\bullet$]
The relevant Fredholm theory can now be established as in~\cite[\S~3]{RobSal95}.

\item[$\bullet$]
The exponential decay of solutions for $s \to \pm \infty$ used in Lemma~\ref{le:compact},
that followed from the Morse Lemma, 
now follows as in~\cite{RobSal01}, see also~\cite[Lemma~2.11]{Sal.lectures}.

\item[$\bullet$]
At the end of Case 1 in the proof of Proposition~\ref{p:DSsur} 
we used unique continuation for ordinary differential equations. 
To establish the unique continuation result needed now one can use a technique of Agmon--Nirenberg as in~\cite[Lemma~3.3]{RobSal01}.
\end{itemize} 


\subsection{A compactness theorem for non-local perturbations of the Cauchy--Riemann operator}
\label{ss:nonlocalcomp}

In this section we prove the compactness result used in the previous paragraph. 
As before we assume that $M = \CC^n$.
Suppose that $J_t(z)$ for $z \in \CC^n$ and $t \in S^1=\RR/\ZZ$ is a smooth family 
of almost complex structures on~$\CC^n$ which is allowed to depend as well smoothly on the variable~$t \in S^1$. 
We again abbreviate the Sobolev spaces
$$
H_\ell \,:=\, W^{\ell,2}(S^1,\CC^n)
$$
and suppose that we are given a vector field $\mathcal{V}$ in the space $\mathcal{C}^k(H,H)$ defined
in~\eqref{e:CkH}, with $k \geqslant 2$.
For example if $V_t$ is a smooth vector field on~$\CC^n$ which may also depend smoothly on $t \in S^1$
we can define such a vector field by 
$$
\mathcal{V}(z)(t) = V_t(z(t)), \quad z \in H_1,\,\,t \in S^1.
$$
However, we do not require that $\mathcal{V}$ be of this form. In particular, $\mathcal{V}$ can be non-local 
in the sense that for a loop $z \in H_1$ and a time $t \in S^1$ the value $\mathcal{V}(z)(t) \in \CC^n$ 
depends on the whole loop~$z$ and not just on the point~$z(t)$ and the time~$t$. 
In fact, that $\mathcal{V}$ is allowed to be non-local is the main novelty of the discussion in this section. 
For $T>0$ abbreviate $I_T = [-T,T]$ and 
$$
\mathcal{C}^k(I_T,H) \,=\,
C^0(I_T,H_k) \cap C^1(I_T,H_{k-1}) \cap \dots \cap C^{k-1}(I_T,H_1) .
$$
We are interested in solutions $w \in \mathcal{C}^k(I_T,H)$
of the following non-locally perturbed Cauchy--Riemann equation on the finite cylinder $I_T \times S^1$:
\begin{equation} \label{cr}
\partial_s w+J_t(w) \partial_t w = \mathcal{V}(w).
\end{equation}
The main result of this section is the following compactness statement.
\begin{theorem} \label{t:compact}
Suppose that $(\mathcal{V}_\nu)$ is a sequence in $\mathcal{C}^k(H,H)$ converging to~$\mathcal{V}$,
and that $(w_\nu) \subset \mathcal{C}^k(I_T,H)$ 
is a sequence of solutions of
\begin{equation} \label{crnu}
\partial_s w_\nu+J_t(w_\nu) \partial_t w_\nu = \mathcal{V_\nu}(w_\nu)
\end{equation}
for which there exists a constant~$C>0$ such that 
$$
\| w_\nu(s)\|_{H_k} \leqslant C \quad \mbox{ for all }\; \nu \in \NN, \; s \in I_T .
$$
Then a subsequence of $(w_\nu)$ converges in $\mathcal{C}^k(I_{T-1},H)$ to a solution~$w$ of~\eqref{cr}.
\end{theorem}

\proof
For notational convenience we assume that~$k=2$. 
We are given a sequence $(w_\nu)_{\nu \geqslant 1} \subset C^0(I_T,H_2) \cap C^1(I_T,H_1)$ 
of solutions of~\eqref{crnu}, and by assumption there exists a constant~$C>0$ such that 
$$
\| w_\nu(s)\|_{H_2 } \leqslant C, \quad
\| \mathcal{V}_\nu(w_\nu(s))\|_{H_2} \leqslant C, \quad
\| D \mathcal{V}_\nu(w_\nu(s)) \|_{\mathcal{L}(H_1,H_1)} \leqslant C, 
\quad \nu \in \NN, \; s \in I_T .
$$
We need to find a subsequence of $(w_\nu)$ that converges in $C^0(I_{T-1},H_2) \cap C^1(I_{T-1},H_1)$ 
to a solution of~\eqref{cr}.
The proof requires the following regularity result.
\begin{proposition} \label{regular}
Suppose that $w \in C^0(I_T,H_2) \cap C^1(I_T,H_1)$ is a solution of~\eqref{cr}. 
Then 
$w \in L^2(I_{T-1},H_3) \cap W^{1,2}(I_{T-1},H_2) \cap W^{2,2}(I_{T-1},H_1)$.
\end{proposition}

\proof
We first explain the proof of the proposition under the simplifying assumption that the almost
complex structure is constant and given by multiplication with~$i$, so that $w$ is a solution of the equation
\begin{equation} \label{crs}
\partial_s w+i \partial_t w = \mathcal{V}(w).
\end{equation}
The proof of this case is inspired by the proof of \cite[Theorem\,3.8]{RobSal95}. 
For $c \in \RR$ let 
$$
A_c \colon H_{k+1} \to H_k, \quad \xi \mapsto i \tfrac{d}{dt} \xi + c \2 \xi.
$$
Note that if $c \notin 2 \pi \ZZ$, then $A_c$ is an isomorphism simultaneously for all $k \in \NN_0$. 
We further abbreviate 
$$
\mathcal{V}_c \in C^0(H_2,H_2) \cap C^1(H_1,H_1), \quad 
\mathcal{V}_c(z) = \mathcal{V}(z)+cz,\,\,z \in H_1.
$$
We can now rewrite~\eqref{crs} as
$$
\partial_s w+A_c w \,=\, \mathcal{V}_c(w).
$$
For
$$
\xi \,:=\, \partial_s w \in C^0(I_T,H_1) \cap C^1(I_T,H_0)
$$
we have
\begin{equation} \label{e:dsxi}
\partial_s \xi+A_c \xi \,=\, D \mathcal{V}_c(w) \xi
\end{equation}
where
$$
D \mathcal{V}_c(z) = D \mathcal{V}(z)+c \cdot \mathrm{id}|_{H_1}, \quad z \in H_1.
$$
Choose a smooth cutoff function $\beta \in C^\infty(I_T,[0,1])$ such that
$$
\beta (s) \,=\,
\left\{\begin{array}{ccl}
1 & \mbox{if} & s \in [-T+1, T-1], \\ [0.2em]
0 & \mbox{if} & s \in [-T,-T+\frac 12] \cup [T-\frac 12,T].
\end{array}\right.
$$
Then
$$
\xi^\beta := \beta \1 \xi \,\in\, C^0(I_T,H_1) \cap C^1(I_T,H_0)
$$
has compact support in~$(-T, T)$. Pick further
$\rho \colon \RR \to [0, \infty)$ such that $\rho(s)=0$ for $|s| \geqslant 1$ and $\int_\RR \rho=1$.
For $\delta>0$ set
$$
\rho_\delta(s) \,=\, \tfrac{1}{\delta} \,\rho \big( \tfrac{s}{\delta} \big).
$$
For $0 < \delta < \frac 12$ we use the notation
$$
\xi^\beta_\delta \,:=\, \rho_\delta * \xi^\beta \in  C^0(I_T,H_1) \cap C^1(I_T,H_0).
$$
Because $\delta < \frac 12$ it holds that $\xi^\beta_\delta$ still has compact support 
in~$(-T,T)$.
For $c \notin 2 \pi \ZZ$ we have
$$
\xi \,=\, -A_c^{-1} \partial_s \xi + A_c^{-1} D \mathcal{V}_c(w) \xi.
$$
Hence we get
\begin{eqnarray} \label{e:xibd}
\xi^\beta_\delta 
&=& 
\rho_\delta* \big(-A_c^{-1} \beta \1 \partial_s \xi+A_c^{-1} D \mathcal{V}_c(w) \xi^\beta \big) \notag \\
&=& 
-(\partial_s \rho_\delta) * (A_c^{-1}\xi^\beta) + \rho_\delta *A_c^{-1} 
                    \Big( (\partial_s \beta) \xi + D_c \mathcal{V}(w) \xi^\beta \Big) . 
\end{eqnarray}
From this formula we deduce that
$$
\xi_\delta^\beta \,\in\, C^0(I_T,H_2) \cap C^1(I_T,H_1).
$$
We introduce the Hilbert spaces
$$
\mathcal{H} \,:=\, L^2(I_T,H_1), \quad 
\mathcal{W} \,:=\, L^2(I_T,H_2) \cap W^{1,2}(I_T,H_1).
$$
The inner product on $\mathcal{W}$ is given by
$\langle \zeta_1, \zeta_2 \rangle_\mathcal{W} =
\langle \zeta_1, \zeta_2 \rangle_{L^2(I_T,H_2)} + \langle \zeta_1, \zeta_2 \rangle_{W^{1,2}(I_T,H_1)}$
where
$$
\langle \zeta_1, \zeta_2 \rangle_{W^{1,2}(I_T,H_1)} \,=\, 
\int_{-T}^T \langle \zeta_1(s), \zeta_2(s) \rangle_{H_1} ds + 
\int_{-T}^T \big\langle \tfrac{d}{dt} \zeta_1(s), \tfrac{d}{dt} \zeta_2(s) \big\rangle_{H_1} ds .
$$
Between these spaces we have the operator
$$
D_c \colon \mathcal{W} \to \mathcal{H}, \quad \zeta \mapsto \partial_s \zeta+A_c \zeta.
$$
Using \eqref{e:xibd} we compute
\begin{eqnarray*}
D_c \1 \xi_\delta^\beta 
&=& (\partial_s \rho_\delta) * \xi^\beta - A_c \Big((\partial_s \rho_\delta)
* (A_c^{-1} \xi^\beta) \Big) + A_c \Big( \rho_\delta * A_c^{-1} \Big( (\partial_s \beta) \xi
  + D_c \mathcal{V}(w) \xi^\beta \Big) \Big) \\
&=& ( \partial_s \rho_\delta) * \xi^\beta - (\partial_s \rho_\delta) * \xi^\beta + \rho_\delta * 
    \Big( (\partial_s \beta) \xi + D_c \mathcal{V}(w) \xi^\beta \Big) \\
&=& \rho_\delta * \Big( (\partial_s \beta) \xi + D_c \mathcal{V}(w) \xi^\beta \Big) .
\end{eqnarray*}
In particular, with the constant 
$\kappa_0 := \left\| (\partial_s \beta) \xi + 
D_c \mathcal{V}(w) \xi^\beta \right\|_{\mathcal{H}}$ 
independent of~$\delta$ we have
\begin{equation} \label{est1}
\| D_c \1 \xi_\delta^\beta \|_{\mathcal{H}} \,\leqslant\, \kappa_0.
\end{equation}
For $\zeta \in \mathcal{W}$ we compute
\begin{eqnarray*}
\|D_c \zeta\|_{\mathcal{H}}^2 
\,=\, \int_{-T}^T \| \partial_s \zeta\|_{H_1}^2 ds + 
                         2 \int_{-T}^T \langle \partial_s \zeta, A_c \zeta \rangle_{H_1} ds
+ \int_{-T}^T \|A_c \zeta\|^2_{H_1} ds .
\end{eqnarray*}
If $\zeta$ has compact support in $(-T,T)$ we obtain, using integration by parts,
$$
\int_{-T}^T \langle \partial_s \zeta, A_c \zeta \rangle_{H_1}ds \,=\,
-\int_{-T}^T \langle \zeta, \partial_s(A_c \zeta) \rangle_{H_1}ds \,=\,
-\int_{-T}^T \langle \zeta, A_c \partial_s \zeta \rangle_{H_1}ds \,=\,
-\int_{-T}^T \langle A_c\zeta, \partial_s \zeta \rangle_{H_1}ds ,
$$
and so
$$
\int_{-T}^T \langle \partial_s \zeta, A_c \zeta \rangle_{H_1}ds \,=\, 0,
$$
whence
\begin{equation}\label{est2}
\| D_c \zeta\|_{\mathcal{H}}^2 \,=\, \int_{-T}^T \| \partial_s \zeta\|_{H_1}^2 ds +
\int_{-T}^T \|A_c \zeta\|^2_{H_1} ds \,\geqslant\, \frac{\|\zeta\|^2_{\mathcal{W}}}{\kappa_1^2}
\end{equation}
for some constant $\kappa_1>0$ depending only on~$c$. Putting \eqref{est1} and \eqref{est2} together we obtain
$$
\| \xi^\beta_\delta\|_{\mathcal{W}} \,\leqslant\, \kappa_1 \1 \kappa_0.
$$
Since the bound $\kappa_1 \1 \kappa_0$ is independent of $\delta$, there exists a sequence $\delta_\nu \to 0$ for
$\nu \to \infty$ such that $\xi^\beta_{\delta_\nu}$ converges weakly in~$\mathcal{W}$
to some $\xi^\beta_0 \in \mathcal{W}$. Since $\xi^\beta_{\delta_\nu}$ converges strongly
in~$\mathcal{H}$ to~$\xi^\beta$, we conclude that $\xi^\beta = \xi^\beta_0 \in \mathcal{W}$. 
Hence
$$
\xi \in L^2(I_{T-1},H_2) \cap W^{1,2}(I_{T-1},H_1) ,
$$
implying that  
$$
w \in W^{1,2}(I_{T-1},H_2) \cap W^{2,2}(I_{T-1},H_1).
$$
Taking advantage of
$$
w=-A_c^{-1} \partial_s w+A_c^{-1} \mathcal{V}(w)
$$
we deduce that 
$$
w \in L^2(I_{T-1},H_3) \cap W^{1,2}(I_{T-1},H_2) \cap W^{2,2}(I_{T-1},H_1).
$$
This proves the proposition in the case where the almost complex structure is constant and given by multiplication 
with~$i$. 

The general case, where the almost complex structure is not constant, follows similarly by a change of the framing. Changing the framing leads to an error term which is of lower order and does not affect the argument.
This was for example used in~\cite{RobSal01}. Here is how this works. 
Suppose that 
$$
\partial_s w + J_t(w) \partial_t w \,=\, \mathcal{V}(w).
$$
Choose a smooth family of matrices such that
$$
\Phi_t(z) J_t(z) \,=\, i \1 \Phi_t(z), \quad z \in \CC^n,\,\,t \in S^1.
$$
Set
\begin{eqnarray*}
\xi(s)(t) &:=& \Phi_t \big( w(s)(t) \big)\partial_s w(s)(t) \\
&=& \Phi_t \big( w(s)(t) \big)
\Big( -J_t \big( w(s)(t) \big) \partial_t w(s)(t) + \mathcal{V} \big( w(s)(t) \big) \Big) \\
&=& -i \1 \Phi_t \big( w(s)(t) \big) \partial_t w(s)(t) + \Phi_t \big( w(s)(t) \big) \mathcal{V} \big( w(s)(t) \big) .
\end{eqnarray*}
We introduce the map
$$
\phi \colon C^0(H_2,H_2) \cap C^1(H_1,H_1) \to C^0(H_2,H_2) \cap C^1(H_1,H_1)
$$
which for a vector field $\mathcal{V} \in C^0(H_2,H_2) \cap C^1(H_1,H_1)$ is given by
$$
\phi(\mathcal{V})(z)(t) = \Phi_t(z(t)) \circ \mathcal{V}(z)(t), \quad z \in H_1,\,\,t \in S^1.
$$
Using that
\begin{eqnarray*}
\Phi_t \big( w(s)(t) \big) \partial_s \partial_t w(s)(t) &=&
\Phi_t \big( w(s)(t) \big) \partial_t \Big( \Phi_t^{-1} \big( w(s)(t) \big) \xi(s)(t) \Big) \\
&=&
-\Big( \dot \Phi_t\big( w(s)(t) \big) + D \Phi_t \big( w(s)(t) \big) \partial_t w(s)(t) \Big) \partial_s w(s)(t) \\
& & + \partial_t \xi(s)(t)
\end{eqnarray*}
we compute
\begin{eqnarray*}
\partial_s \xi(s)(t) &=& -i \Big( D \Phi_t \big(w(s)(t) \big) \partial_s w(s)(t) \Big) \partial_t w(s)(t)
-i \Phi_t \big( w(s)(t) \big) \partial_s \partial_t w(s)(t) \\
& & + D (\phi(\mathcal{V})) \partial_s w(s)(t) \\
&=& -i \Big( D \Phi_t \big( w(s)(t) \big) \partial_s w(s)(t) \Big) \partial_t w(s)(t) \\
& & +i \Big( \dot \Phi_t\big( w(s)(t) \big) + D \Phi_t \big( w(s)(t) \big) \partial_t w(s)(t) \Big) \partial_s w(s)(t) 
-i \partial_t \xi(s)(t) \\
& & + D( \phi( \mathcal{V}) ) \partial_s w(s)(t) .
\end{eqnarray*}
Writing $\Psi_c(w) \in C^0(I_T,H_1)$ for
\begin{eqnarray*}
\Psi_c(w)(s)(t) &=& -i\Big(D \Phi_t\big(w(s)(t)\big)\partial_s w(s)(t)\Big) \partial_t w(s)(t) \\
& & \Big( i \1 \dot \Phi_t + i D \Phi_t \,\partial_t w(s)(t) + D(\phi(\mathcal{V}))  \Big) \big( w(s)(t) \big)
\, \partial_s w(s)(t)  + c \2 \xi(s)(t)
\end{eqnarray*} 
we then get the compact expression
$$
\partial_s \xi + A_c \xi \,=\, \Psi_c(w) ,
$$
similar to~\eqref{e:dsxi}.
Using this, the proof for the general case now proceeds like the one in the special case where~$J$ was constant. 
\proofend

The following lemma is an elaboration of \cite[Lemma\,3.6]{RobSal95}.
\begin{lemma} \label{complem}
For $T>0$ the inclusion 
\begin{eqnarray*}
\iota \colon L^2(I_T,H_3) \cap W^{1,2}(I_T,H_2) \cap W^{2,2}(I_T,H_1) \,\to\, 
                     C^0(I_T,H_2) \cap C^1(I_T,H_1) 
\end{eqnarray*}
is a compact operator. 
\end{lemma}

\proof
For $N \in \NN$ abbreviate by $V_N \subset H_0 = L^2(S^1,\CC^n)$
the subspace of finite Fourier series of the form
$$
z = \sum_{k=-N}^N z_k \,e^{2 \pi i k}, \quad z_k \in \CC^n.
$$
Note that $V_N$ is finite-dimensional. Indeed, its real dimension is $2n(2N+1)$. Let
$$
\pi_N \colon H_0 \to V_N
$$
be the orthogonal projection. Note that the Fourier basis is a common orthogonal basis of the Sobolev spaces
$H_k = W^{k,2}(S^1,\CC^n)$ for every $k \in \NN_0$. In particular, the restriction 
$$
\pi_N|_{H_k} \colon H_k \to V_N
$$
coincides with the orthogonal projection of $H_k$ to~$V_N$. Let
$$
\Pi_N \colon  L^2(I_T,H_3) \cap W^{1,2}(I_T,H_2) \cap
W^{2,2}(I_T,H_1) \to W^{2,2}(I_T,V_N), \quad w \mapsto \pi_N \circ w.
$$
Because $V_N$ is finite-dimensional, the inclusion
$$
I_N \colon W^{2,2}(I_T,V_N) \to C^1(I_T,V_N)
$$
is a compact operator. We finally abbreviate by
$$
J_N \colon C^1(I_T,V_N) \to C^0(I_T,H_2) \cap C^1(I_T,H_1)
$$
the inclusion. Denote by
$$
\iota_N \colon L^2(I_T,H_3) \cap W^{1,2}(I_T,H_2) \cap W^{2,2}(I_T,H_1) \to 
C^0(I_T,H_2) \cap C^1(I_T,H_1)
$$
the composition of these three maps,
$$
\iota_N \,:=\, J_N \circ I_N \circ \Pi_N.
$$
Because $I_N$ is compact and the other two maps are continuous, $\iota_N$ is a compact operator. 
We show that $\iota_N$ converges to~$\iota$ in the norm topology as $N \to \infty$. 
To see this we argue by contradiction and assume that there exists a constant~$c>0$ such that 
for every $N \in \NN$ there exists $w_N \in L^2(I_T,H_3) \cap W^{1,2}(I_T,H_2) \cap W^{2,2}(I_T,H_1)$ 
with the property that 
\begin{equation} \label{max1}
\max \Big\{ \| (\iota-\iota_N) w_N \|_{C^0(I_T,H_2)},\, \| (\iota-\iota_N) w_N \|_{C^1(I_T,H_1)} \Big\} \,=\, 1
\end{equation}
but
\begin{equation} \label{max2}
\max \Big\{ \| w_N \|_{L^2(I_T,H_3)},\, \|w_N\|_{W^{1,2}(I_T,H_2)},\, \|w_N\|_{W^{2,2}(I_T,H_1)} \Big\} \,\leqslant\, c.
\end{equation}
Choosing $c$ larger if necessary we can assume that $\frac{1}{4c^2} \leqslant T$.
From~\eqref{max1} we deduce that 
\begin{equation} \label{fall1}
\| (\iota-\iota_N) w_N\|_{C^0(I_T,H_2)} = 1
\end{equation}
or
\begin{equation} \label{fall2}
\| (\iota-\iota_N)w_N \|_{C^1(I_T,H_1)} = 1.
\end{equation}
We first discuss case \eqref{fall1}. In this case there exists $s \in I_T$ such that 
$$
\| (\mathrm{id}-\pi_N) w_N(s) \|_{H_2} = 1.
$$
For $s' \in I_T$ with $|s'-s| \leqslant \frac{1}{4c^2}$ we have
$$
|s'-s| \,\leqslant\, \frac{1}{4c^2} \,\leqslant\, \frac{1}{4 \|w_N\|^2_{W^{1,2}(I_T,H_2)}} \,\leqslant\,
\frac{1}{4 \|(\mathrm{id}-\pi_N)w_N \|^2_{W^{1,2}(I_T,H_2)}} .
$$
We can thus estimate
\begin{eqnarray*}
\| (\mathrm{id}-\pi_N) w_N(s') \|_{H_2} &\geqslant&  
\| (\mathrm{id}-\pi_N) w_N(s) \|_{H_2} -
\bigg| \int_s^{s'} \|(\mathrm{id}-\pi_N) \partial_\sigma w_N (\sigma) \|_{H_2} d\sigma \bigg| \\
&\geqslant&
1- \| (\mathrm{id}-\pi_N)w_N \|^2_{W^{1,2} (I_T, H_2)} \sqrt{|s'-s|} \\
&\geqslant& 1-\tfrac{1}{2} \\
&=& \tfrac{1}{2}.
\end{eqnarray*}
In particular,
$$
\| (\mathrm{id}-\pi_N)w_N(s')\|_{H_3} \geqslant \pi N
$$
from which we deduce that
$$
\| w_N \|_{L^2(I_T,H_3)} \,\geqslant\, \sqrt{\frac{\pi^2 N^2}{4 c^2}} \,=\, \frac{\pi N}{2 c}
$$
which contradicts \eqref{max2} as soon as $N > \tfrac{2 c^2}{\pi}$. 
This contradiction shows that case~\eqref{fall1} cannot occur, 
and case~\eqref{fall2} leads to a contradiction in a similar way. 
This shows that $\iota_N$ converges to $\iota$ as $N \to \infty$. In particular, 
$\iota$ arises as the limit of compact operators and is therefore compact itself. 
This finishes the proof of the lemma. 
\proofend

\ni
{\it Proof of Theorem~\ref{t:compact}, for $k=2$.} 
Under the hypothesis of Theorem~\ref{t:compact}, and in view of equation~\eqref{crnu},
the proof of Proposition~\ref{regular} reveals that there exists a bounded subset 
$$
\mathcal{U} \,\subset\,  L^2(I_{T-1},H_3) \cap W^{1,2}(I_{T-1},H_2) \cap W^{2,2}(I_{T-1},H_1)
$$
such that 
$w_\nu|_{I_{T-1}} \in \mathcal{U}$ for all $\nu \in \NN$.
Hence by Lemma~\ref{complem} a subsequence of $w_\nu|_{I_{T-1}}$ converges in
$C^0(I_{T-1},H_2) \cap C^1(I_{T-1},H_1)$ to~$w$, that clearly solves~\eqref{cr}. 
\proofend


\end{document}